\documentclass[reqno, psamsfonts, 11pt]{amsart}
\usepackage{mathabx}    
\usepackage{mathrsfs}  
\usepackage[backref,colorlinks]{hyperref} 
\usepackage{amsmath,amssymb,amsthm}
\usepackage[abbrev,msc-links]{amsrefs}   
\usepackage[utf8]{inputenc}
\usepackage{tikz}
\usepackage{amsmath}
\usepackage{enumerate}
\usepackage{setspace}
\usepackage{lineno}
\usepackage{float}
\theoremstyle{plain}
\newtheorem{thm}{Theorem}[section]

\newtheorem{prop}[thm]{Proposition}

\newtheorem*{clm*}{Claim}
\newtheorem{claim}[thm]{Claim}

\theoremstyle{definition}

\newtheorem{dfn}[thm]{Definition}

\numberwithin{equation}{section}

\normalsize                              
\let\theta=\vartheta
\let\rho=\varrho
\let\phi=\varphi

\def\cT{{\mathcal T}}

\def\od{\overline{d}}

\usepackage{enumitem} 

\DeclareMathOperator{\ex}{ex}

\let\polishlcross=\l
\def\l{\ifmmode\ell\else\polishlcross\fi}

\def\qand{\quad\text{and}\quad}

\def\eps{\varepsilon}

\def\tU{\tilde U}
\def\tZ{\tilde Z}

\theoremstyle{plain}
\newtheorem{theorem}{Theorem} 
\let\epsilon=\varepsilon

\newtheoremstyle{note}
  {4pt}
  {4pt}
  {\sl}
  {}
  {\itshape}
  {.}
  {.5em}
  {}

\begin{document}


\title[]{Tur\'an number of complete multipartite graphs in multipartite graphs}
\thanks{
Jie Han is partially supported by a Simons Collaboration Grant 630884 and National Natural Science Foundation of China (12371341).
Yi Zhao is partially supported by NSF grant DMS 2300346 and Simons Collaboration Grant 710094.}
%
\author{Jie Han}
\address{School of Mathematics and Statistics and Center for Applied Mathematics\\ Beijing Institute of Technology\\ Beijing\\ China}
\email{han.jie@bit.edu.cn}
\author{Yi Zhao}
%
%
\address
{Department of Mathematics and Statistics, Georgia State University, Atlanta, GA 30303}
\email{yzhao6@gsu.edu}

\begin{abstract}
In this paper we study a multi-partite version of the Erd\H{o}s--Stone theorem.
Given integers $r<k$ and $t\ge 1$, let $\ex_k(n, K_{r+1}(t))$ be the maximum number of edges of $K_{r+1}(t)$-free $k$-partite graphs with $n$ vertices in each part, where $K_{r+1}(t)$ is the complete $(r+1)$-partite graph with $t$ vertices in each part.
We determine the exact value of $\ex_k(n, K_{r+1}(t))$ for $t\le 3$, $r<k\le 2r$ and sufficiently large $n$.
{We also characterize all extremal graphs for $r, k$ such that $r$ divides $k$, analogous to a result of Erd\H os and Simonovits on forbidding $K_{r+1}(t)$ in general graphs.}
\end{abstract} 

\maketitle


\section{Introduction}
Generalizing Mantel's theorem from 1907 \cite{Mantel}, Tur\'an's theorem from 1941 \cite{MR0018405} started the systemetic study of Extremal Graph Theory. Given a graph $F$, let $\ex(n, F)$ denote the largest number of edges in a graph not containing $F$ as a subgraph (called $F$-free). Let $K_r$ denote the complete graph on $r$ vertices and $T_{r}(n)$ denote the complete $r$-partite graph on $n$ vertices with $\lfloor n/r \rfloor$ or $\lceil n/r \rceil$ in each part (known as the Tur\'an graph); and $t_r(n)$ be the size of $T_{r}(n)$.
Tur\'an's theorem \cite{MR0018405} states that $\ex(n, K_{r+1})= t_r(n)$ for all $n\ge r\ge 1$ and in addition, $T_{r}(n)$ is the unique \emph{extremal graph}. 

Let $K_{t_1, \dots, t_r}$ denote the complete $r$-partite graph with parts of size $t_1, \dots, t_r$ and write
$K_r(t) = K_{t, \dots, t}$ with $r$ parts. 
A celebrated result of Erd\H os and Stone \cite{MR0018807} determines $\ex(n, K_{r+1}(t))$ asymptotically: 
\[
\ex(n, K_{r+1}(t)) = t_r(n) + o(n^2) = \left(1 - \frac1r \right)\frac{n^2}2 + o(n^2).
\]
Erd\H os \cite{MR0232703} and Simonovits \cite{MR0233735} independently improved the error term above to $O(n^{2- 1/t})$.
Simonovits \cite{MR0233735} also showed that any extremal graph for $K_{r+1}(t)$ can be obtained from $T_r(n)$ by adding or removing $O(n^{2 - 1/t})$ edges. 
Later Erd\H os and Simonovits \cite{MR0292707} determined the structure of extremal graphs for $K_{r+1}(t)$ for $t\le 3$ as follows. 

\begin{theorem}\cite{MR0292707}
\label{thm:ES}
For $t\le 3$, every extremal graph $G$ for $K_{r+1}(t)$ has a vertex partition $U_1, \dots, U_r$ such that
\begin{itemize}
\item $G[U_i, U_j]$ is complete for all $i\ne j$,
\item $G[U_i] = n/r + o(n)$,
\item $G[U_1]$ is extremal for $K_{t, t}$, and
\item $G[U_2], \dots, G[U_r]$ are extremal for $K_{1, t}$.
\end{itemize}
\end{theorem}
The restriction $t\le 3$ in Theorem~\ref{thm:ES} comes from our knowledge on $\ex(n, K_{t, t})$. A well-known open problem in Extremal Graph Theory is proving $\ex(n, K_{t, t}) = \Omega(n^{2 - 1/t})$ and this is only known for $t\le 3$.

{Extremal problems with multipartite graphs as host graphs have been studied since 1951, when 
Zarankiewicz proposed the study of the largest number of edges in a $K_{s, t}$-free bipartite graph.}
Let $\mathcal G(n_1, \dots, n_k)$ denote the family of $k$-partite graphs with $n_1, \dots, n_k$ vertices in its parts and write $\mathcal G_k(n)= \mathcal G(n, \dots, n)$ with $k$ parts. Given a graph $F$, define 
$\ex(n_1, \dots, n_k; F)$ as the largest number of edges in $F$-free graphs from $\mathcal{G}(n_1, \dots, n_k)$, and let 
$\ex_k(n, F)= \ex(n, \dots, n; F)$ (with $k$ parts).
(Trivially $\ex_k(n, F)= \binom k2 n^2$ if the chromatic number $\chi(F)>k$.) 
In 1975 Bollob\'{a}s, Erd\H{o}s, and Szemer\'{e}di \cite{MR0389639} investigated several Tur\'an-type problem for multipartite graphs. 
Applying a simple counting argument, they showed that, for any $n, k, r\in \mathbb{N}$ with $k> r$,  
\begin{align}
    \label{eq:exkn}
    \ex_k(n, K_{r+1})= t_r(k) n^2
\end{align}
 
 The main results of \cite{MR0389639} were on the minimum degree that forces a copy of $K_r$ in the graphs of $\mathcal{G}_k(n)$.  This problem has been intensively studied and resolved when $k=r$ (frequently in its complementary form concerning \emph{independent transversals}) \cites{MR2195582, MR1860440, MR1208805, MR2246152}, and more recently for $k>r$ \cite{MR4496024} (certain cases are still unsolved). There are several other extremal results for multipartite graphs. 
Bollob\'{a}s, Erd\H{o}s, and Straus \cite{MR0379256} determined $\ex(n_1, \dots, n_r; K_r)$ for all $n_1, \dots, n_r$ and $r$. Let $kK_r$ denote $k$ vertex-disjoint copies of $K_r$. The problem $\ex(n_1, \dots, n_r; kK_r)$ were studied more recently \cites{MR3724116,MR4425157,MR4595892}
and settled by 
Chen, Lu, and Yuan~\cite{CLY23} when $n_1, \dots, n_r$ are sufficiently large ($k, r$ are arbitrary but fixed). The minimum pair density of multipartite graphs that forces a clique was 
studied by Bondy, Shen, Thomass\'e, and Thomassen \cite{bondy} and Pfender \cite{pfender}.

\bigskip
In this paper we study $\ex_k(n, K_{r+1}(t))$, the multi-partite version of the Erd\H os--Stone theorem and Theorem~\ref{thm:ES}.
To give the precise value of $\ex_k(n, K_{r+1}(t))$, we need the following definition.
Given $a, t, n_1,\dots, n_a\in \mathbb N$,
let $z_t^{(a)}(n_1,\dots, n_a)$ be the $a$-partite Zarankiewicz number for $K_{t,t}$, that is, the maximum number of edges in a $K_{t,t}$-free bipartite graph with part sizes $n_1,\dots, n_a$.
For simplicity, we write $z_t^{(a)}(n)$ if $n_1=\cdots=n_a$.
We also write $z_t(m,n)$ for $z_t^{(2)}(m, n)$.

Theorem~\ref{thm:ES} says that any extremal graph for $K_{r+1}(t)$ is the \emph{join} of an extremal graph for $K_{t,t}$ and $r-1$ extremal graphs for $K_{1,t}$.
Inspired by this, a natural guess of an extremal graph for its $r$-partite analogue is as follows.
Suppose $k=ar+b$ with $0\le b<r$. We start with an $n$-blowup of $T_r(k)$, which has $r$ classes and each class has either $a$ or $a+1$ parts. We add a $K_{t, t}$-free graph to one class with the most number of parts and $\{K_{1, t}, K_{2, 2}\}$-free graphs to the remaining $r-1$ classes.
If $r < k\le 2r$, then it is easy to see (as in the proof of Theorem~\ref{thm:lblb}) that this graph is $K_{r+1}(t)$-free and has 
$t_r(k)n^2+z_t(n,n)+(k-r-1)(t-1)n$
edges, and therefore 
\[
\ex_k(n, K_{r+1}(t))\ge t_r(k)n^2+z_t(n,n)+(k-r-1)(t-1)n.
\]

In this paper we 
first improve this lower bound when $t\ge 2$ and $r<k\le 2r$.
\begin{theorem}
\label{thm:lblb}
Suppose $t\ge 2$, $r<k\le 2r$ and $n\ge 8t^2$.
Then $\ex_k(n, K_{r+1}(t))\ge g(n,r,k,t)$, where
\[
g(n,r,k,t) := t_r(k)n^2 + z_t(n, n) + (t-1)(k-r-1)n + \min\{k-r-1, 2r-k\}\left\lfloor \frac{(t-1)^2}4 \right\rfloor.
\]
\end{theorem}
The following is our main result, in which we prove a matching upper bound when $t\in \{2,3\}$ and $n$ is sufficiently large.
\begin{theorem}
\label{thm:lb}
For any $t\in \{2,3\}$ and integers $r$, there exists $n_0=n_0(r)\in \mathbb N$ such that for $n\ge n_0$ and $r<k\le 2r$,
we have $\ex_k(n, K_{r+1}(t))=g(n,r,k,t)$.
\end{theorem}
In fact, we conjecture that $\ex_k(n, K_{r+1}(t))=g(n,r,k,t)$ holds for all $t\ge 2$, $r<k\le 2r$, and sufficiently large $n$.


A natural question is whether our result can be extended to larger values of $t$ and $k$.
For larger value of $t$, although we can use $z_{t}(n,n)$ without knowing its precise value, we need several properties of this function in our proof. K\H{o}v\'ari, S\'os, Tur\'an \cite{MR0065617} showed that $z_{t}(n,n) =  O(n^{2-1/t})$ for $t\ge 2$ and proving a matching lower bound is a well-known open problem:
\begin{enumerate}[label=($Z$)]
\item $z_{t}(n,n) =  \Omega(n^{2-1/t})$ for $t\ge 2$. 
\label{item:zest}
\end{enumerate}
It was shown \cites{MR0200182,MR0223262} that \ref{item:zest} holds for $t=2, 3$.
In addition, we will need the following properties.
\begin{enumerate}[label=($E\arabic*$)]
\item For any $a\in \mathbb{N}$, there exists $\delta>0$ such that for large $n$, $z_{t}^{(a+1)}(n) - z_{t}^{(a)}(n) \ge {\delta n^{2-1/t}}. $\label{item:za=1gap}
\item for any $\eps\in (0,1]$, there exists $\delta>0$ such that for large $n$,  \label{item:zagap}
\[
z_{t}(n,n) - z_{t}((1-\eps)n, n) \ge {\delta n^{2-1/t}}.
\] 
\item $z_{t}(m, n) - z_{t}(m-1, n) \ge t-1$. 
\label{item:zgap}
\end{enumerate}
We can easily verify \ref{item:za=1gap}--\ref{item:zgap} for $t=2, 3$. First, a proof of \ref{item:za=1gap} from \ref{item:zest} for all $t$ is given in Section~\ref{sec:33}. Second, 
when $t=2, 3$, \ref{item:zagap} follows from $z_t(m,n) \le (1+o(1))mn^{1-1/t}$ by F\"uredi \cite{MR1395691} and $z_t(n,n)\ge (1-o(1))n^{2-1/t}$ for $t=2, 3$.
Third,~\ref{item:zgap} holds trivially because adding a vertex with $t-1$ edges to a $K_{t, t}$-free graph will not create a copy of $K_{t, t}$.

We do not know whether similar properties hold when $k > 2r$: in \ref{item:zagap} we must deal with the $\lceil k/r\rceil$-partite Zarankiewicz number; we also need to replace $t-1$ by $\Omega(a^2t)$ in \ref{item:zgap}, which seems out of reach at present.



Theorems~\ref{thm:lblb} and~\ref{thm:lb} show that our problem is more complex than its non-partite counterpart, Theorem~\ref{thm:ES}.
Finally, we show that, when $r$ divides $k$, this additional complexity does not exist, and we give an analogue of Theorem~\ref{thm:ES}, modulo the existence of a set of constantly many exceptional vertices.

\begin{theorem}\label{thm7}
For $r, k\in \mathbb N$ with $r\mid k$ and $t=2,3$, there exist $C_0, n_0\in \mathbb N$ such that the following holds for $n\ge n_0$.
Let $G$ be a $K_{r+1}(t)$-free $k$-partite graph with $n$ vertices in each part and $\ex_k(n, K_{r+1}(t))$ edges.
Then there is a vertex partition of $G$ into $U_1,\dots, U_r$, each consisting of exactly $k/r$ parts of $G$, and a vertex set $Z\subseteq V(G)$ with $|Z|\le C_0$ such that
\begin{itemize}
\item $G[U_i\setminus Z, U_j\setminus Z]$ is almost complete for all $i\ne j$,
\item $G[U_1\setminus Z]$ is $K_{t, t}$-free, and
\item $G[U_2\setminus Z], \dots, G[U_r\setminus Z]$ are $K_{1, t}$-free.
\end{itemize}
\end{theorem}

{Showing $Z=\emptyset$ in Theorem~\ref{thm7} requires an $a$-partite analogue of~\ref{item:zagap}, which is unknown for $a\ge 3$.}

For the rest of this paper we only consider $r\ge 2$, as it is easy to see that Theorems~\ref{thm:lblb}--\ref{thm7} hold for $r=1$.


\medskip
\noindent

\medskip
\noindent
\textbf{Notation.}
Given integers $n\ge m\ge 1$, let $[n]=\{1, \dots, n\}$ and $[m, n]=\{m, m+1, \dots, n\}$.
We omit floors and ceilings unless they are crucial, e.g., we may choose a set of $\eps n$ vertices even if our assumption does not guarantee that $\eps n$ is an integer.

When $X, Y\subseteq V(G)$ intersect, $E_G(X, Y)$ is defined as the collection of ordered pairs in $(x,y)\in X\times Y$ such that $\{x, y\}\in E(G)$.
Write $e_G(X, Y) = |E_G(X, Y)|$.
For a vertex $v$ in $G$, let $N(v, X)=N(v)\cap X$ and $d(v, X)=|N(v, X)|$.
Moreover, given $X\subseteq V(G)$, let $e(X, G)$ be the number of edges of $G$ incident to the vertices of $X$.
Given two graphs $G$ and $H$ on a common vertex set $V$, $G\cap H$ denotes a graph on $V$ with $E(G\cap H)=E(G)\cap E(H)$.
Given a $k$-partition $\{V_1,V_2,\dots,V_k\}$, a set $S$ is called \emph{crossing} if $|S\cap V_i|\le 1$, $i\in [k]$.

When we choose constants $x,y>0$, $x\ll y$ means that for any $y> 0$ there
exists $x_0> 0$ such that for any $x< x_0$ the subsequent statement holds.
Hierarchies of other lengths are defined similarly.
Furthermore, all constants in the hierarchy are positive and for a constant appearing in the form $1/s$, we always mean to choose $s$ as an integer.

\section{Proof of Theorem~\ref{thm:lblb}}

In this section we prove Theorem~\ref{thm:lblb}, that is, $\ex_k(n, K_{r+1}(t))\ge g(n,r,k,t)$ for $r<k\le 2r$.
{Our proof needs a $t$-regular $K_{2,2}$-free bipartite graph with $n$ vertices in each part. 
It is well known (see \cite{MR101250}) that such graph exists for infinitely many $n\in \mathbb{N}$ with $n\ge t^2$. 
The following proposition from \cite{MR2519934}*{Section 2}
allows $n$ to be any integer that is at least $8t^2$.} 

\begin{prop}\cite{MR2519934}
\label{fact:cons1}
For $t\ge 1$ and $n\ge 8t^2$,
there exists a $t$-regular $K_{2,2}$-free bipartite graph with $n$ vertices in each part.
\end{prop}

 Now we prove our lower bound on $\text{ex}_k(n, K_{r+1}(t))$ stated in Theorem~\ref{thm:lblb}.

\begin{proof}[Proof of Theorem~\ref{thm:lblb}]
First assume $r<k\le 2r$ and $t\ge 2$.
Let $V_1,\dots, V_{k}$ be disjoint sets of size $n$ and, if $k<2r$, let $V_{k+1},\dots, V_{2r}$ be empty sets.
Let $G'$ be a complete $r$-partite graph with parts $V_1\cup V_{r+1}, V_2\cup V_{r+2}, \dots, V_{r}\cup V_{2r}$.
Moreover, we add to $G'$ a maximum $K_{t,t}$-free bipartite graph with bipartition $V_1\cup V_{r+1}$, and a $(t-1)$-regular $K_{2,2}$-free bipartite graph on $V_{i}\cup V_{i+r}$ for $2\le i\le k-r$ (the existence of such graph is guaranteed by Proposition~\ref{fact:cons1}).
The resulting graph is $K_{r+1}(t)$-free because a copy of $K_{r+1}(t)$ has at most $2t-1$ vertices in $V_1\cup V_{r+1}$ and at most $t$ vertices in $V_{i}\cup V_{i+r}$ for $2\le i\le k-r$.
This graph has $t_r(k)n^2 + z_{t}(n) + (t-1)(k-r-1)n$ edges, and thus, $\ex_k(n, K_{r+1}(t)) \ge t_r(k)n^2 + z_{t}(n) + (t-1)(k-r-1)n$. This proves the theorem when $t=2$ or $k=2r$.

Now assume $r<k< 2r$ and $t\ge 2$. Let $b=k-r$. Our goal is to give a better construction that shows
\[
 \ex_k(n, K_{r+1}(t)) \ge t_r(k)n^2 + z_t(n, n) + (t-1)(b-1)n + \min\{b-1, r-b\}\left\lfloor \frac{(t-1)^2}4 \right\rfloor.
\]
Let $V_{i,j}$, $(i,j)\in [r]\times [2]$ be vertex sets, where $V_{b+1, 2}, \dots, V_{r, 2}$ are empty sets and other sets have size $n$. Let $G= K(V_{1,1}\cup V_{1, 2}, \dots, V_{r, 1}\cup V_{r, 2})$ be the complete $r$-partite graph with parts $V_{1,1}\cup V_{1, 2}, \dots, V_{r, 1}\cup V_{r, 2}$. Thus $e(G)= t_r(k)n^2$. 

We first revise the partition as follows.
Let $t':=\lceil (t-1)/2\rceil$ and $b':=\min\{b-1, r-b\}$. 
Let $\{V_{i,j}', (i,j)\in [r]\times [2]\}$ be obtained from $\bigcup V_{i,j}$ by moving a set $S_{i,1}$ of $t'$ vertices from $V_{i, 1}$ to $V_{i+b-1, 1}$, and moving a set $S_{i,2}$ of $t'$ vertices from $V_{i, 2}$ to $V_{i+b-1, 2}$, for every $i\in [2, b'+1]$.
For $i\in [r]$, let $U_i:=V_{i, 1}' \cup V_{i,2}'$ and $H := K(U_1,\dots, U_{r})\cap K(V_{1,1},\dots, V_{r,1}, V_{1,2}, \dots, V_{r,2})$. 
Let $H'$ be obtained from $H$ by adding
\begin{itemize}
\item a $K_{t,t}$-free bipartite graph on $U_1$ of size $z_t(n,n)$, and
\item a maximum $\{K_{1,t}, K_{2,2}\}$-free bipartite graph on $U_i$ for $i\in [2,b]$.
\item $2t'$ vertex-disjoint copies of $K_{1,t-1}$ on each $U_i$, $i\in [b+1, b+b']$, where the centers of stars are the $2t'$ vertices of $S_{i-b+1, 1}\cup S_{i-b+1, 2}$; then add a copy of $K_{t', t'}$ on $S_{i-b+1, 1}\cup S_{i-b+1, 2}$. 
See Figure \ref{fig:enter-label}.
\end{itemize}

\begin{figure}
    \centering
    \includegraphics[width=0.85\linewidth]{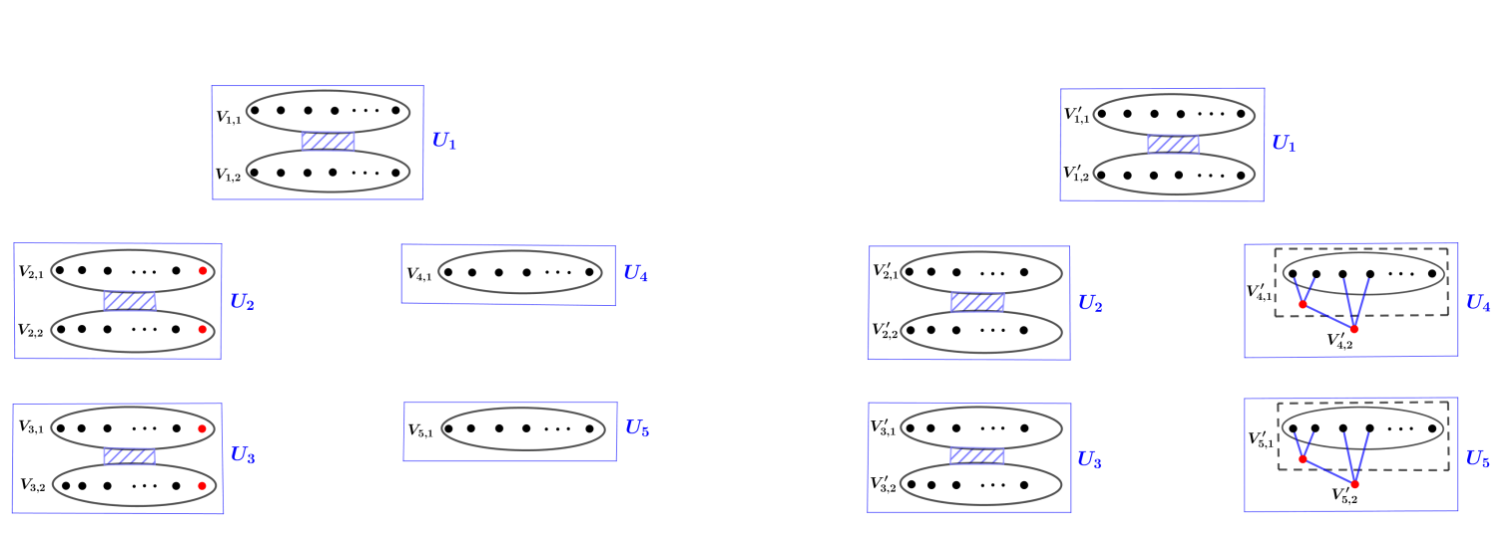}
    \caption{The lower bound construction for $K_{r+1}(t)$ with $r=5$, $k=8$ and $t=3$. The left figure is the standard construction similar to the one given in Theorem~\ref{thm:ES}; the right figure is the construction presented in our proof of Theorem~\ref{thm:lb}.}
    \label{fig:enter-label}
\end{figure}

For $i\in [2, b]$, Proposition~\ref{fact:cons1} implies that each $H'[U_i]$ is $(t-1)$-regular and thus for $i\in [2, b'+1]$, $H'[U_i]$ has $(n-t')(t-1)$ edges, and for $i\in [b'+2, b]$ it has $n(t-1)$ edges; for $i\in [b+1, b+b']$, each $H'[U_i]$ has $2t'(t-1)+(t')^2$ edges.
Therefore, the number of edges in $H'\setminus H$ is
\begin{align*}
    & z_{t}(n,n)+(n-t')(t-1)b'+n(t-1)(b-b'-1)+b'(2t'(t-1)+(t')^2)\\
    =& z_t(n,n)+(t-1)(b-1)n+b'(t')^2+b't'(t-1)\\
    =& z_t(n,n)+(t-1)(b-1)n+2b'(t')^2+b' t'(t-1-t')\\
    =& z_t(n,n)+(t-1)(b-1)n+2b'(t')^2+ \left\lfloor \frac{(t-1)^2}4 \right\rfloor.
\end{align*}

We claim that $H$ contains $t_r(k)n^2 - 2b'(t')^2$ edges.
Indeed, for every $i\in [2, b'+1]$,  the vertices of $S_{i,1}$ moved from $V_{i, 1}$ to $V_{i+b-1, 1}$ lose $n$ edges to $V_{i+b-1, 1}$ and gain $n-t'$ edges to $V_{i, 2}$, thus having a net loss $(t')^2$ edges between $U_i$ and $U_{i+b-1}$; the same holds for $S_{i,2}$. 
Thus, our claim holds after summing over all $b'$ such rows, which implies that 
$e(H') = t_r(k)n^2 + z_{t}(n,n) + (t-1)(b-1)n + b'\lfloor \frac{(t-1)^2}4 \rfloor$. 


At last, we show that $H'$ is $K_{r+1}(t)$-free.
Recall that by construction, every $U_i$ is triangle-free, $U_1$ is $K_{t,t}$-free, $U_2,\dots, U_b$ are $\{K_{1,t}, K_{2,2}\}$-free, and each of $U_{b+1},\dots, U_{U_{b+b'}}$ induces vertex-disjoint copies of $K_{1,t-1}$ whose centers are joined by copies of $K_{t',t'}$.
Suppose $K$ is a copy of $K_{r+1}(t)$ in $H'$. 
By construction, $K$ contains at most $2t-1$ vertices in $U_1$.
We claim that $|V(K)\cap (U_i\cup U_{i+b-1})|\le 2t$ for $i\in [2,b'+1]$ and $|V(K)\cap U_j|\le t$ for $i\in [b'+2,b]\cup [b+b'+1,r]$, which will lead to a contradiction with $|V(K)|=(r+1)t$.
Indeed, for $i\in [2,r]$, if $K$ contains at least $t+1$ vertices in $U_i$, then these vertices induce either a copy of $K_{1,t}$ or $K_{2,2}$.
By construction, this is only possible for $i\in [b+1,b+b']$ and in that case $V(K)\cap U_i$ must intersect both $S_{i-b+1, 1}$ and $S_{i-b+1,2}$.
Furthermore, there exists $v\in V(K)\cap S_{i-b+1, 1}$ and $v'\in V(K)\cap S_{i-b+1, 2}$ such that $v$ and $v'$ are in different color classes of $K$. Since $v$ and $v'$ have no common neighbor in $U_i\cup U_{i-b+1}$, $K$ induces at most two classes on $U_i\cup U_{i-b+1}$.
Thus, $|V(K)\cap (U_i\cup U_{i-b+1})|\le 2t$ and we are done.
\end{proof}

\noindent\textbf{Remark.}
Note that for $t=2$, both our constructions give the same value, which means that the extremal graph is not unique.
Indeed, it is easy to see that one can construct $b'+1$ different ones -- as we can move vertices for a subset of the $b'$ rows.

\section{Stability and proof outline}

Let us first consider extremal graph for $\ex_k(n, K_{r+1})$. 
Given $r, k\in \mathbb N$ with $k> r$, write $k=ar+b$ for $0\le b\le r-1$.
By Tur\'an's theorem, the Tur\'an graph $T_r(k):=K_{a,\dots, a, a+1,\dots, a+1}$ (with $b$ parts of size $a+1$ and $r-b$ parts of size $a$) is the unique largest $K_{r+1}$-free graph on $k$ vertices. 
The following definition shows that there are many extremal graphs for $\ex_k(n, K_{r+1})$. 

\begin{dfn}
\label{dfn1}
Let $\cT_{r,k}(n)$ be the collection of $k$-partite graphs with parts $V_1, \dots, V_k$ of size $n$ defined as follows.
If $b>0$, we arbitrarily divide $V_{ar+1}, \dots, V_k$ into $r$ sets $W_1,\dots, W_{r}$ (some of them may be empty) such that each $W_i$ is a subset of $V_j$ for some $j$;
if $b=0$, then let $W_1,\dots, W_r$ be empty sets.
Now let $T$ be the $r$-partite graph with parts $U_1, \dots, U_r$ such that 
\[
U_i=W_i\cup Z_i, \text{ where } Z_i:=V_{(i-1)a+1}\cup \cdots\cup V_{ia}, 
\]
obtained from the complete $r$-partite graph $K(U_1,\dots, U_{r})$ by removing edges between $W_i$ and $W_{i'}$, $i\ne i'$, whenever $W_i, W_{i'}\subseteq V_j$ for some $j$
(in other words, $T=K(U_1,\dots, U_{r})\cap K(V_1,\dots, V_k)$).  
\end{dfn}

Since $T$ is $r$-partite, it is $K_{r+1}$-free.
Let $U=\bigcup_{i\in [r]}U_i$ and $W=\bigcup_{i\in [r]}W_i$.
Note that $e_T(U)=\binom r2 a^2n^2$ while the number of edges of $T$ incident to $W$ is equal to $|W|(k-a-1)n=b(k-a-1)n^2$.
Since $t_r(k)=\binom r2 a^2+b(k-a-1)$, it follows that $e(T)=t_r(k) n^2$. By \eqref{eq:exkn}, $T$ is an extremal graph for $K_{r+1}$.\footnote{Indeed, Chen, Lu, and Yuan~\cite{CLY23} showed that $\cT_{r,k}(n)$ are \emph{all} the extremal graphs for $\ex_k(n, K_{r+1})$ and described all extremal graphs for $\ex(n_1, \dots, n_k; K_{r+1})$.} 


\subsection{A stability theorem}
We need the following stability result for $\ex_k(n, K_{r+1}(t))$.
Given two graphs $G, H\in \mathcal G_k(n)$ on the same parts $V_1, \dots, V_k$, we say that $G$ and $H$ are $\gamma$-close if $|E(G)\triangle E(H)|\le \gamma n^2$.

\begin{theorem}\label{thm:sta2}
For any $k, r, t\in \mathbb N$ and any $\gamma >0$, there exist $\eps>0$ and $n_0\in \mathbb N$ such that the following holds for every integer $n\ge n_0$.
Suppose $G\in \mathcal G_k(n)$ is $K_{r+1}(t)$-free and $e(G)\ge (t_r(k)-\eps)n^2$.
Then $G$ is $\gamma$-close to a member of $\cT_{r,k}(n)$.
In particular, we have $e(G)\le (t_r(k)+\gamma)n^2$.
\end{theorem}

In the earlier version of this paper, we gave a self-contained proof of Theorems~\ref{thm:sta2}. Here we derive it from a stronger result of Chen, Lu and Yuan~\cite{CLY23}*{Theorem 1.5}, in which they provide more structural information.
Their definition and result are more general by allowing the parts $V_i$ of $G$ to have different sizes. Here we only state their definition and result that we need for the balanced case.

\begin{dfn}
[Stable partition]
Let $k\ge r\ge 2$ be integers.
Let $\mathcal P=\{P_1,\dots, P_{r}\}$ and $\mathcal V=\{V_1,\dots, V_k\}$ where $|V_i|=n$ for all $i\in [k]$ be partitions of a common vertex set.
For $i\in [r],j\in [k]$, a set $W=P_i\cap V_j$ is called an integral part if $W=V_j$, and called a partial part otherwise.
We say that $\mathcal P$ is stable to $\mathcal V$, if each of $P_1,\dots, P_r$ has the same number of integral parts and at most one partial part.
\end{dfn}

\begin{dfn}
[$(X,\varepsilon)$-stable]
Let $n,r,k$ be integers with $n\ge k\ge r\ge 2$.
For a given spanning subgraph $G$ of $K_{k}(n)$, let $\mathcal P=\{P_1,\dots, P_{r}\}$ be a partition of $V(G)$ and $\mathcal V=\{V_1,\dots, V_k\}$ be the natural $k$-partition.
Given  $\varepsilon\in (0,1)$ and a set $X\subseteq V(G)$ of size at most $\varepsilon n$,  we say that $\mathcal P$ is an $(X,\varepsilon)$-stable partition if
\begin{itemize}
    \item $G-X$ is $\varepsilon$-close to $K(P_1,\dots, P_r)-X$,
    \item $\{P_1\setminus X,\dots, P_{r}\setminus X\}$ is stable to $\{V_1\setminus X,\dots, V_{k}\setminus X\}$.
\end{itemize}
\end{dfn}

\begin{theorem}\cite{CLY23}*{Theorem 1.5}
\label{thm:CLY}
Let $F$ be a graph with chromatic number $r+1\ge 3$.
For every $\varepsilon>0$ there exists $\delta>0$ and integer $n_0>0$ such that the following holds for $n\ge n_0$.
Let $G$ be an $F$-free subgraph of $K_{k}(n)$ with $k>r$ such that $e(G) \ge \ex_k(n, F) - \delta n^2$.
Then, $G$ has an $(X, \varepsilon)$-stable partition $\{P_1,\dots, P_{t-1}\}$ {for some set $X\subseteq V(G)$ of size at most $\varepsilon n$}.
\end{theorem}

Now Theorem~\ref{thm:sta2} follows from Theorem~\ref{thm:CLY} by setting $F=K_{r+1}(t)$ and noticing that i) a partition is stable if and only if it is a partition required in Definition~\ref{dfn1} and ii) as $|X|\le \varepsilon n$, $X$ is incident to at most $\varepsilon n^2$ edges that we have no control on.

We remark that although Theorem~\ref{thm:CLY} is stronger,  the additional structural information is not immediately useful to us.
For example, the set $X$ introduced in Theorem~\ref{thm:CLY} is the set of atypical vertices of $G$ (see also the proof outline in next section).
However, we need a stronger control and indeed we need to distinguish two kinds of atypical vertices.
Identifying them from $X$ is almost equivalent to identifying them from $V(G)$. So we choose to present and use the easier and more classical version, Theorem~\ref{thm:sta2}, in this paper.





{
\subsection{Outline of the proofs.}
Now we give an outline of our proofs.
Let $G\in \mathcal G_k(n)$ be $K_{r+1}(t)$-free and has the maximum number of edges.
Since $e(G)>t_r(k)n^2$, we can assume that $G$ is $\gamma$-close to some $T\in \cT_{r,k}(n)$.
Since $e(G)$ is maximum, we can easily derive a minimum degree condition by symmetrization arguments.

Next we define \emph{atypical} vertices.
Roughly speaking, there are two types of \emph{atypical} vertices: the first type of vertices, denoted by $Z''\cup W''$, are the ``wrong'' ones that do not exist in $\cT_{r,k}(n)$; the second type of vertices, 
denoted by $(W'_i\setminus W_i) \cup \bigcup_{j\ne i} Z_j^i$ for $i\in [r]$,
are the vertices that are not in $U_i$ but \emph{behave} like the vertices of $U_i$, in other words, they are in the wrong place.
We temporarily ignore the first type of atypical vertices because there are only a constant number of them (see~\ref{item:Wsizes} and~\ref{item:Zsizes}) and they contribute only $O(n)$ to $e(G)$.
For the second type of atypical vertices, there are only $o(n)$ of them (see~\ref{item:Wsizes} and~\ref{item:Zsizes}) and we move them to appropriate rows and redefine our partition as $\tU_1, \dots, \tU_r$ (see~\eqref{eq:partition1}).
A key observation is that $Z_j^i\neq \emptyset$ (namely, there is a vertex in $Z_j$ but behaves as a vertex in $Z_i$) only if $|W_j| \ge (1-o(1))n$.
}

{
Now we estimate $e(G)$. 
We split $E(G)$ into $E_G(\tU_1), \dots, E_G(\tU_r)$, and $E(G')$, where $G':=G\cap K(\tU_1,\dots, \tU_{r})$.
We have a relatively good estimate of $e(G')$ (see Claim~\ref{clm:eG'}) taking into account that the partition is no longer balanced.
In contrast, due to the second type of atypical vertices, we can only show that each $G[\tU_i]$ is ``almost'' $K_{t,t}$-free (see Claim~\ref{clm:Ktt}).
Similarly, we show that all but at most one rows are ``almost'' $K_{1,t}$-free (see Claim~\ref{clm:K1t}).
Assuming that $e_G(\tU_1)$ is the largest among all $e_G(\tU_i)$, $i\in [r]$, we can use these properties and give an upper bound of $e(G)$.
Next we show that $\tU_1$ has no atypical vertices (Claim~\ref{clm:clean}), and thus $e_G(\tU_1)\le z_t^{(a+1)}(n)$.
We further refine our estimate on $\tU_i$, $i>1$, and show that each second type atypical vertex contributes at most a constant number of edges to $E(G)\setminus E(G')$ (Claim~\ref{clm:clean0}).
In summary, $\tU_1$ is indeed $K_{t,t}$-free, $e_G(\tU_i)=O(n)$ for $i>1$, and $|Z''\cup W''|=O(1)$, from which we conclude the proof of Theorem~\ref{thm7}.
}

{
To prove Theorem~\ref{thm:lb}, we refine earlier estimates as follows.
We first show that $|W_1'| = (1-o(1))n$ and $Z''\cup W''=\emptyset$, where we use~\ref{item:zagap}.
The rest of the proofs are further refinements of our estimates.
In particular, we show that if $|W_i'|$ is not too small, then $\tU_i$ essentially contains no vertex from other rows.
}

\subsection{Two quick proofs}
\label{sec:33}
We first derive \ref{item:za=1gap} from~\ref{item:zest}.

\begin{proof}
[Proof of \ref{item:zest} $\Rightarrow$ \ref{item:za=1gap}]
Given $a+1$ sets $V_1, \dots, V_{a+1}$ of size $n$, we define an $(a+1)$-partite graph $G$ on $V_1,  \dots, V_{a+1}$ as follows.
Let $V_2'$ be a set of $n$ vertices consisting of $\lfloor n/2\rfloor$ vertices from $V_1$ and $\lceil n/2\rceil$ vertices from $V_2$.
We place an extremal graph $G'$ for $z_{t}^{(a)}(n)$ on $V_2', V_3, \dots, V_{a+1}$, in other words, $G'$ is an $a$-partite $K_{t,t}$-free graph with $z_{t}^{(a)}(n)$ edges.
Next we add a maximum bipartite $K_{t,t}$-free graph $G''$ on the remaining vertices of $V_1$ and $V_2$.
By \ref{item:zest}, $e(G'')\ge z_t^{(2)}(\lfloor n/2\rfloor) \ge \delta n^{2-1/t}$ for some $\delta >0$. 
Thus $G=G'\cup G''$ is $K_{t,t}$-free and $e(G) = e(G') + e(G'') \ge z_{t}^{(a)}(n) + \delta n^{2-1/t}$.
This gives~\ref{item:za=1gap}. 
\end{proof}

We need the following simple proposition. 
\begin{prop}
\label{lem:badvtx}
Given $r, t\in \mathbb N$ and reals $\gamma, \eps>0$ such that $\eps^2 > 3r^2t^2\gamma$, and let $n$ be sufficiently large.
Suppose $G$ is a $K_{r+1}(t)$-free graph with vertex partition $V=U_1\cup\cdots\cup U_{r}$ such that $|U_i|\ge n$ for $i\in [r]$ and $d(U_i, U_j)\ge 1-\gamma$, $i, j\in [r]$, $i\neq j$.
Let $X\subseteq V$ be the set of vertices $v$ satisfies that $d(v, U_i) \ge \eps |U_i|$ for all $i\in [r]$.
Then $|X|\le 2(t-1)\eps^{-r t}$.
\end{prop}

\begin{proof}
We call a copy of $K_{r}{(t)}$ in $G$ \emph{useful} if it consists of exactly $t$ vertices from each of $U_1, \dots, U_{r}$.
We first show that for every $v\in X$, $N(v)$ contains many useful copies of $K_{r}{(t)}$.
Indeed, $d(U_i, U_j)\ge 1-\gamma$ for every $i, j\in [r]$, $i\neq j$ implies that $G[U_i, U_j]$ has at most $\gamma |U_i||U_j|$ non-edges.
Since $d(v, U_i) \ge \eps |U_i|$ for all $i\in [r]$, take $W_i\subseteq N(v)\cap U_i$ of size exactly $\eps |U_i|$.
We can find $\prod_{i\in [r]}\binom{\eps |U_i|}{t}$ $rt$-sets which consists of $t$ vertices from each $U_i$, amongst which, at most
\[
\sum_{i, j\in [r], i\neq j}\gamma |U_i||U_j| \cdot \left(\prod_{i'\in [r]}\binom{\eps |U_{i'}|}{t}\right)\cdot \frac{t}{\eps |U_i|} \frac{t}{\eps |U_j|} \le \frac{r^2t^2\gamma}{\eps^2} \prod_{i'\in [r]}\binom{\eps |U_{i'}|}{t}
\]
of them contain crossing non-edges.
Therefore, $N(v)$ contains at least
\[
\left(1-\frac{r^2t^2\gamma}{\eps^2}\right) \prod_{i'\in [r]}\binom{\eps |U_{i'}|}{t} \ge \frac{\eps^{rt}}2 \prod_{i'\in [r]}\binom{|U_{i'}|}{t}
\]
useful copies of $K_{r}{(t)}$, where we used that $r^2t^2\gamma \eps^{-2} < 1/3$.
Since $G$ is $K_{r+1}(t)$-free, each useful copy $K$ of $K_{r}{(t)}$ is in $N(v)$ for at most $t-1$ choices of $v\in X$.
Double counting on the number of pairs $(v, K)$ such that $K\subseteq N(v)$ is useful, we obtain that 
\[
|X| \frac{\eps^{rt}}2 \prod_{i'\in [r]}\binom{|U_{i'}|}{t} \le (t-1) \prod_{i'\in [r]}\binom{|U_{i'}|}{t},
\]
which gives $|X|\le 2(t-1)\eps^{-rt}$.
\end{proof}

\section{Main Proofs}

Given integers $1\le s\le t$ and sufficiently large $m, n$, 
K\H{o}v\'ari, S\'os, Tur\'an \cite{MR0065617} showed that $z(m,n,s,t)\le C m n^{1-1/s}$ for some $C=C(t)>0$, that is, a bipartite graph $G$ with parts of size $m$ and $n$ has at most $C m n^{1-1/s}$ edges if $G$ has no copy of $K_{s,t}$ where the part of size $s$ is in the part of $G$ of size $m$. 
%
{We start the proof with the general setting, that is, for $k=ar+b$ with $0\le b < r$. After we conclude the proof of Theorem~\ref{thm7}, we focus on the case $k\le 2r$.}

\begin{proof}
[Proofs of Theorems~\ref{thm:lb} and~\ref{thm7}]
Suppose $t=2,3$ and thus \ref{item:zest} holds, that is, $z_t^{(2)}(n) \ge c n^{2-1/t}$ for some $c >0$. 
Take $C=C(t)$ as in the K\H{o}v\'ari--S\'os--Tur\'an result in the previous paragraph.
We choose constants 
\[
1/n \ll \gamma \ll \eps \ll \eps' \ll 1/k, 1/t, c, C.
\]

Suppose $G$ is $K_{r+1}(t)$-free and has the maximum number of edges, that is, $e(G)=\ex_k(n, K_{r+1}(t))$.
Suppose further that $e(G)> g(n, r, k,t) > t_r(k) n^2$.
By Theorem~\ref{thm:sta2}, $G$ is $\gamma$-close to some $T\in \cT_{r,k}(n)$ such that $T=K(U_1,\dots, U_{r})\cap K(V_1,\dots, V_k)$, where $U_1,\dots, U_{r}$ is a partition of $V(T)$ satisfying the following properties:
\begin{itemize}
    \item $U_i=W_i\cup Z_i$ such that $Z_i=V_{(i-1)a+1}\cup \cdots\cup V_{ia}$, and
    \item $W_i=\emptyset$ if $b=0$ and $W_i$ is a subset of $V_{q_i}$ for some $q_i$ with $ar < q_i \le k$ otherwise.
\end{itemize}
For simplicity, we write $z_t(n) = z_t^{(\lceil k/r \rceil)}(n)$.

\medskip
The fact that $G$ is $\gamma$-close to $T$ gives the following observation.
\begin{enumerate}[label=(D0)]
\item for any $i\in [r]$, there exists $B_i\subseteq U_i$ of size at most $2\sqrt\gamma n$ such that for any $v\in U_i\setminus B_i$ and $A\subseteq \bigcup_{j\in [r]\setminus \{i\}}U_j$ satisfying that none of the vertices of $A$ is in the same cluster as $v$ is, we have $\overline d(v, A)\le \sqrt\gamma n$. \label{item:Deg0}
\end{enumerate}
To see it, fix $i\in [r]$ and write $U^*:=\bigcup_{j\ne i} U_j$. Since $G$ is $\gamma$-close to $T$, we have
\[
e_G(Z_i, U^*) \ge |Z_i| |U^*| - \gamma n^2, \quad \text{and} \quad e_G \left(W_i, U^* \setminus V_{q_i} \right) \ge |W_i||U^*\setminus V_{q_i}| - \gamma n^2.
\]
Let $B_i'\subseteq Z_i$ be the set of vertices $v$ such that $\overline d(v, U^*) >\sqrt\gamma n$, and
$B_i''\subseteq W_i$ be the set of vertices $w$ such that $\overline d(w, U^*\setminus V_{q_i}) > \sqrt\gamma n$.
The displayed line above implies that $|B_i'|\le \sqrt\gamma n$ and $|B_i''|\le \sqrt\gamma n$.
Now~\ref{item:Deg0} holds by setting $B_i=B_i'\cup B_i''$.

\subsection*{Minimum degree}
For $i\in [k]$, let $N_i:=N_T(u_i)$ for some $u_i\in V_i$.
Note that this is well-defined as the vertices of $V_i$ share the same neighborhood in $T$.
Using the maximality of $e(G)$, we derive that for every $u\in V_i$, $i\in [k]$
\[
     d_G(u) \ge d_T(u) - 2t\gamma n.
\]
Indeed, since $G$ is $\gamma$-close to $T$, that is, $|E(G)\triangle E(T)|\le \gamma n^2$, for each $i\in [k]$, we can greedily pick distinct $u_1,\dots, u_t\in V_i$, such that $|N_G(u_j)\triangle N_i|\le \gamma n^2/(n-j) \le 2\gamma n$, for $j\in [t]$.
Let $N_i':=\bigcap_{j\in [t]}N_G(u_j)$ and note that $|N_i'\triangle N_i|\le 2t\gamma n$.
In particular, $|N_i'|\ge |N_i| - 2t\gamma n$.
Now for a contradiction suppose there is $u\in V_i$ such that $d_G(u) < d_T(u) - 2t\gamma n = |N_i| - 2t\gamma n$.
Then we replace $N_G(u)$ by $N_i'$, that is, we disconnect all the edges of $u$ in $G$ and connect $u$ to the vertices of $N_i'$.
Thus, we obtain a $k$-partite graph on the same vertex set as $G$ and has more edges than $G$.
Therefore, by the maximality of $G$, this new graph contains a copy of $K_{r+1}(t)$, denoted by $K$.
Clearly, $K$ must contain the vertex $u$, as $G$ is $K_{r+1}(t)$-free.
Moreover, $K$ must miss at least one vertex from $u_1,\dots, u_t$, say $u_j$, because the set $\{u, u_1,\dots, u_t\}$ is independent in $G$ and $K$ has independence number $t$.
However, as the neighborhood of $u$ $N_i'$ is a subset of $N_G(u_j)$, we can replace $u$ by $u_j$ and still get a copy of $K_{r+1}(t)$, which is in $G$, a contradiction.

Therefore, comparing with the degrees in $T$, we derive that for any vertex $u$,
\begin{equation}
\label{eq:mindeg}
d_G(u) \ge \begin{cases} 
(k-a)n - |W_i| - 2t\gamma n, \,\,\,\,\hfill \text{ if $u\in Z_i$ for $i\in [r]$,} \\
(k-1-a)n - 2t\gamma n, \,\,\,\,\hfill \text{ if } u\in W. 
\end{cases}
\end{equation}

\subsection*{Atypical vertices}
In this step we identify a set of atypical vertices, that is, those behave differently from the majority of the vertices.
Let $W:=\bigcup_{i\in [r]}W_i = V_{ar+1}\cup \cdots\cup V_k$.
We define $W'':=\{v\in W: d(v, Z_j)\ge \eps n, \text{ for all }j\in [r]\}$ and $W_i':=\{v\in W: d(v, Z_i) < \eps n\}$. 
Then we have $W= W''\cup W_1'\cup \cdots \cup W_{r}'$.
Next, for $i\in [r]$, let $Z'':=\bigcup_{i\in [r]}Z_i''$, where
\[
Z_i'':=\{v\in Z_i: d(v, Z_j)\ge \eps n, \text{ for all }j\in [r]\setminus \{i\} \qand d(v, U_i)\ge \eps n\}.
\]
Furthermore, let $Z_i':=Z_i\setminus Z_i''$ 
and write $Z_i'$ as $\bigcup_{j\in [r]} Z_i^j$, where $Z_i^j$, $j\neq i$, consists of the vertices $v\in Z_i$ such that $d(v, Z_j) < \eps n$, and $Z_i^i$ consists of the vertices $v$ such that $d(v, U_i) < \eps n$.
The following are some useful properties of these sets.


\begin{claim}
The following properties hold for all $i\in [r]$.
\begin{enumerate} [label=$(P\arabic*)$]
\item $|W_i'\setminus W_i|\le 2\gamma n$ and $|W''|\le C_0:=2t\eps^{-rt}$.
\label{item:Wsizes}
\item $W= W''\cup W_1'\cup \cdots \cup W_{r}'$ is a partition of $W$.
\item $|Z_i''|\le C_0$, $|Z_i^j|\le \sqrt{\gamma} n$ for $j\ne i$, and $|Z_i^i|\ge (1-\sqrt\gamma)a n$. \label{item:Zsizes}
\item $\bigcup_{j\in [r]} Z_i^j$ is a partition of $Z_i'$. \label{item:Zpartition}
\end{enumerate}
\end{claim}

\begin{proof} 
Recall the definition of $W''$ and that $d(Z_i, Z_j)\ge 1-\gamma$ for distinct $i,j\in [r]$.
Applying Proposition~\ref{lem:badvtx} to the graph $G[W''\cup Z]$ with vertex partition $(U_1, \dots, U_{r})$, we obtain that $|W''|\le C_0:=2t\eps^{-rt}$.
We next show that $|W_i'\setminus W_i|\le 2\gamma n$ for each $i\in [r]$.
Indeed, because $G$ is $\gamma$-close to $T$, we have $e_G(Z_i, W_i'\setminus W_i)\ge an |W_i'\setminus W_i| - \gamma n^2$.
On the other hand, by definition, $e_G(Z_i, W_i'\setminus W_i) < |W_i'\setminus W_i| \cdot \eps n$.
Thus, we get $|W_i'\setminus W_i| < \gamma n /(a-\eps) < 2\gamma n$, verifying~\ref{item:Wsizes}.



To see $(P2)$, suppose there is a vertex $v\in W_i'\cap W_j'$. By definition, $d(v)\le (k-1)n - 2(a-\eps) n < (k-1-a) n- \sqrt\gamma n$, contradicting~\eqref{eq:mindeg}.

Next we show~\ref{item:Zsizes}.  Fix $i\in [r]$.
Since $G$ is $\gamma$-close to $T$, we have $d(Z_j, Z_{j'}) \ge 1-\gamma$ and $d(U_i, Z_j) \ge 1-\gamma$ for distinct $j, j'\in [r]\setminus \{i\}$.
Thus, we can apply Proposition~\ref{lem:badvtx} on $G[U_i\cup \bigcup_{j\neq i} Z_j]$ (with the obvious $r$-partition) and obtain $|Z_i''|\le C_0$.
Moreover, for $i\neq j$, from $d(Z_i, Z_j) \ge 1-\gamma$ we infer $|Z_i^j|\le (\gamma/\eps) n\le \sqrt{\gamma} n$, as $\gamma \ll \eps$.
Therefore, we also get $|Z_i^i|\ge |Z_i| - |Z_i''| - \sum_{j\neq i}|Z_i^j|\ge an - C_0 - (r-1)\gamma n /\eps \ge  (1-\sqrt\gamma)a n$.

Now we show~\ref{item:Zpartition}.
By definition, if $v\in Z_i^i$, then $d(v, U_i) < \eps n$; if $v\in Z_i^j$ for $j\neq i$, then $d(v, Z_j) < \eps n$. 
Thus, we have $Z_i'\subseteq \bigcup_{j\in [r]} Z_i^j$ by definition.
A vertex $v\in Z_i^i \cap Z_i^{j}$, $j\ne i$, satisfies that $d(v)< kn- (|U_i|-\eps n) - (a-\eps)n \le (k - a)n - |W_i| - (1-2\eps)n$, contradicting \eqref{eq:mindeg}.   
A vertex $v\in Z_i^j\cap Z_i^{j'}$ for distinct $j, j'\in [r]\setminus \{i\}$ satisfies that $d(v)< (k-1)n -2(a-\eps)n\le (k-a-2)n+2\eps n$, contradicting \eqref{eq:mindeg} as well. Thus, $\bigcup_{j\in [r]} Z_i^j$ is a partition of $Z_i'$.
%
\end{proof}

For $i\in [r]$, our refined partition is defined by
\begin{equation}
\label{eq:partition1}
\tU_i := \tZ_i \cup W_i' , \text{ where } \tZ_i:= \bigcup_{j\in [r]} Z_j^i .
\end{equation}
Then $V(G)=Z''\cup W''\cup \bigcup_{i\in [r]}\tU_i$.
Note that for any $v\in \tU_i$, we have $d(v, Z_i^i)\le d(v, Z_i)\le \eps n$, and thus $d(v, \tZ_i)\le \eps n + (r-1)\sqrt\gamma n$ by~\ref{item:Zsizes}.

For every $i\in [r]$, note that~\ref{item:Wsizes} implies that $|W_i\setminus W_i'|\le C_0+(r-1)2\gamma n \le 2r\gamma n$, and similarly~\ref{item:Zsizes} implies that $|Z_i\setminus \tZ_i|\le C_0+(r-1)\sqrt\gamma n \le r\sqrt{\gamma}n$.

We now derive a more handy minimum degree condition. For convenience, define $\overline d(v, A)= |A| - d(v, A)$. For $v\in Z_i^i$, we have $\od(v, \tU_i)\ge \od(v, U_i) - |U_i\setminus \tU_i|$.
Since $\od(v, U_i)> an+|W_i|- \eps n$ and $|U_i\setminus \tU_i|\le |Z_i\setminus \tZ_i| + |W_i\setminus W'_i|\le \eps n/2$, we have $\od(v, \tU_i)\ge an+|W_i|-\eps n - \eps n/2$. By \eqref{eq:mindeg}, $\od(v)\le an + |W_i| + \sqrt{\gamma}n$. It follows that $\overline d(v, V\setminus \tU_i )\le 2\eps n$. Now consider $v\in \tU_i\setminus Z^i_i$. The definition of $\tU_i$ implies that $d(v, Z_i)< \eps n$ and $\od(v, Z_i)> an- \eps n$. 
Assume $v\in V_j$. Then $V_j\cap Z_i = \emptyset$ and trivially $\od(v, V_j)=n$. It follows that $\od(v, Z_i\cup V_j) > (a+1)n - \eps n$. 
Hence $\od(v, \tZ_i\cup V_j)\ge \od(v, Z_i\cup V_j) - |Z_i\setminus \tZ_i|> (a+1)n - \frac32\eps n$. On the other hand, either case of \eqref{eq:mindeg} implies that $\od(v)\le (a+1)n + \sqrt{\gamma}n$. Consequently, 
$\overline d(v, V\setminus (\tZ_i\cup V_j))\le 2\eps n$.
In summary, for $i\in [r]$ and $j\in [k]$, 

\smallskip
\begin{enumerate}[label=(Deg)]
\item 
If $v\in Z_i^i$, then $\overline d(v, V\setminus \tU_i)\le 2\eps n$; 
if $v\in (\tU_i \setminus Z^i_i) \cap V_j$, then $\overline d(v, V\setminus (\tZ_i\cup V_j))\le 2\eps n$. \label{item:DegZ}
\end{enumerate}

Next we prove further properties on $Z_i^j$ and $\tZ_j$.

\begin{claim}
\label{clm:Wnew}
If $Z_i^j\neq \emptyset$ for some $i\ne j$, then the following holds.
\begin{enumerate} [label=$(Q\arabic*)$]
\item For $v\in Z_i^j$ and $A\subseteq V(G)\setminus (Z_i\cup Z_j)$, we have $d(v, A)\ge |A| - \eps n - \sqrt\gamma n$.
\label{item:Zijvdeg}
\item $|W_i|\ge (1-\eps-\sqrt\gamma)n$. \label{item:WiLarge}
\item If $|\tZ_j\setminus Z_j|\ge t$, then $|W_j|\le  2t\eps n$. \label{item:Wismall}
\end{enumerate}
\end{claim}

\begin{proof}
Note that $d(v, Z_j)\le \eps n$ and $d(v, Z_i)\le (a-1)n$, that is, $v$ has at least $n+(an-\eps n)=(a+1)n - \eps n$ non-neighbors in $Z_i\cup Z_j$.
On the other hand,~\eqref{eq:mindeg} says that $v$ has at most $an+|W_i|+\sqrt\gamma n$ non-neighbors in $G$.
Combining these two we get that $v$ has at most $|W_i|-n+\eps n + \sqrt\gamma n\le \eps n + \sqrt\gamma n$ non-neighbors outside $Z_i\cup Z_j$, and thus~\ref{item:Zijvdeg} holds.
The fact that $|W_i|-n+\eps n + \sqrt\gamma n\ge 0$ implies~\ref{item:WiLarge}.

For~\ref{item:Wismall}, suppose to the contrary, $|\tZ_j\setminus Z_j|\ge t$ and $|W_j|>  2t\eps n$.
By~\ref{item:Zijvdeg} with $A=W_j$, arbitrary $t$ vertices in $\tZ_j\setminus Z_j$ have at least $|W_j| - t(\eps+\sqrt{\gamma})n \ge t$ common neighbors in $W_j$. We thus obtain a copy of $K_{t,t}$ with one part in $\tZ_j\setminus Z_j$ and the other part in $W_j$ -- denote its vertex set by $B$. 
For any $i'\in [r]\setminus \{j\}$ such that $B\cap Z_{i'}^j\neq \emptyset$, we have $|W_{i'}|\ge (1-\eps-\sqrt\gamma)n$ by~\ref{item:WiLarge}. Since $|W_j|>  2t\eps n$, $W_{i'}$ and $W_j$ do not belong to the same cluster, and thus no vertex of $B$ is in the same cluster that contains $W_{i'}$, which implies that the vertices of $B$ have at least $|W_{i'}| - 2t (2\eps n)\ge n/2$ common neighbors in $W_{i'}$ by~\ref{item:DegZ}.
For any $i''\in [r]\setminus \{j\}$ such that $B\cap Z_{i''}^j= \emptyset$ (and thus $B\cap Z_{i''}= \emptyset$), by~\ref{item:DegZ} we have that the vertices of $B$ have at least $n/2$ common neighbors in $Z_{i''}$.
Because $G$ is $\gamma$-close to $T$, these common neighborhoods, each of size at least $n/2$, have densities close to one between each pair, and thus contain a copy of $K_{r-1}{(t)}$. 
Together with $B$, they form a copy of $K_{r+1}(t)$ in $G$, a contradiction.
\end{proof}

In particular, when $b=0$ (and thus $W_i=\emptyset$ for all $i$), \ref{item:WiLarge} implies that  $Z_i^j=\emptyset$ whenever $i\neq j$. Consequently, 
\begin{equation}
\label{eq:Ui0}
\tU_i=Z_i^i = Z_i\setminus Z'' \quad \text{for all $i\in [r]$ when $b=0$}.
\end{equation}

Let $L\subseteq [r]$ be the set of indices $i$ such that $|W_i|\ge (1-\eps-\sqrt\gamma)n$. \ref{item:WiLarge} and \ref{item:Wismall} imply that
\begin{itemize}
\item for $i\in [r]\setminus L$, we have $Z_i^j= \emptyset$ for $j\ne i$.
\item for $i\in L$, $|\tZ_i\setminus Z_i|\le t-1$ and thus $|\tZ_i|\le an+t-1$.
\end{itemize}

\subsection*{First Estimate on $e(G)$}
Let $G'=G\cap K(\tU_1,\dots, \tU_{r})$. We have $e(G) = e(G') + \sum_{i=1}^r e_G(\tU_i) + e(Z''\cup W'', G)$.
Since $G'$ is $r$-partite, it is $K_{r+1}$-free.
As $G'$ is a subgraph of $G\in \mathcal{G}_k(n)$,
we have $e(G')\le t_r(k)n^2$ (but this is not good enough when $b>0$).
Below we give an upper bound for $e(G')$, which will be used throughout the proof. Recall that $T=K(V_1,\dots, V_{k})\cap K(U_1,\dots, U_{r})$ has precisely $t_r(k) n^2$ edges.

\begin{claim}\label{clm:eG'}
We have $e(G')\le t_r(k)n^2 + \sum_{i\in [r]} (\beta_i - \alpha_i)$, where
\[
\begin{aligned}
\beta_i:=&\sum_{j\in L\setminus \{i\}}|Z_j^i| \left(|\tZ_j\setminus Z_j|+|W_j'|-n + |Z_i\setminus \tZ_i|\right) \text{ and }\\
\alpha_i:=&|\tZ_i\setminus Z_i| |W_{i}'| +e_T(W_i')+e_T(\tZ_i\setminus Z_i).
\end{aligned}
\]
\end{claim}

\begin{proof}
We first obtain $G^{(0)}:=K(Z_1\cup W_1', \dots, Z_{r}\cup W_{r}')\cap K(V_1, \dots, V_k)$ from $T$. During this process, we lose the edges of $T$ between $W_i$ and $W_j$, $j\ne i$, if both ends of the edges are placed in $W'_i$. Thus 
\begin{align}
    \label{eq:eG0}
    e(G^{(0)}) = t_r(k)n^2 - \sum_{i\in [r]} e_T(W_i').
\end{align}

We imagine a dynamic process of obtaining $G'$ from $G^{(0)}$ by recursively moving vertices. To estimate $e(G')$, we track the changes of the edges with respect to complete $r$-partite graphs (but also respecting the $k$-partition of $G$). More precisely, for $l>0$, let
\[
G^{(l)}:=K(Z_1^{(l)}\cup W_1', \dots, Z_{r}^{(l)}\cup W_{r}')\cap K(V_1, \dots, V_k)
\] 
such that the $r$-partition of $G^{(l)}$ can be obtained by moving exactly one vertex from the partition of $G^{(l-1)}$. 
The process terminates after $m:= \sum_{i\in [r]}|\tZ_i\setminus Z_i|$ steps and thus $G'$ is a subgraph of $G^{(m)}$.
Furthermore, throughout the process, we only move vertices from the color classes in $L$ to other color classes.
Therefore, we can give a linear ordering to the members of $L$, and for $i\in L$ we move vertices from $Z_i$ only after we have moved the vertices in color classes $j$ prior to $i$ (denoted by $j<_Li$).
Now, in the $l$-th step, suppose we move $v$ from $Z_{j}^{(l-1)}$ to $Z_{i}^{(l-1)}$, namely, $v\in Z_{j}^{i}$, then the change is
\[
e(G^{(l)}) - e(G^{(l-1)}) = |Z_{j}^{(l-1)}\setminus V_p|+|W_{j}'| - |\tZ_{i}^{(l-1)}| - |W_{i}'|,
\]
where $V_p\ni v$ and $\tZ_{i}^{(l-1)}= Z_{i}^{(l-1)}\setminus V_p$.

Note that we have $|Z_{j}^{(l-1)}\setminus V_p|\le (a-1)n+|\tZ_j\setminus Z_j|$.
Moreover for any $j' <_L j$, we have $Z_{j'}^i\subseteq Z_{i}^{(l-1)}$.
Therefore, we have $|\tZ_{i}^{(l-1)}|\ge an - |Z_i\setminus \tZ_i|+\sum_{j' <_L j} |Z_{j'}^i|$. 
Putting all these together, we get 
\[
e(G^{(l)}) - e(G^{(l-1)}) \le |\tZ_j\setminus Z_j|+|W_j'|-n + |Z_i\setminus \tZ_i|-\sum_{j' <_L j} |Z_{j'}^i| - |W_{i}'|.
\]
Recalling that we moved $v$ from $Z_{j}^{(l-1)}$ to $Z_{i}^{(l-1)}$ at the $l$-th step, we obtain
\[
e(G') - e(G^{(0)})\le \sum_{l=1}^{m} \left( |\tZ_j\setminus Z_j|+|W_j'|-n + |Z_i\setminus \tZ_i|-\sum_{j' <_L j} |Z_{j'}^i| - |W_{i}'| \right),
\]
where $i, j$ depends on $l$.
Since $m= \sum_{i\in [r]}|\tZ_i\setminus Z_i|$, we have 
\[
\begin{aligned}
&\sum_{l=1}^{m}(|\tZ_j\setminus Z_j|+|W_j'|-n+ |Z_i\setminus \tZ_i| - |W_{i}'|) \\
 = \, &\sum_{i\in [r]}\sum_{j\in L\setminus \{i\}}|Z_j^i| (|\tZ_j\setminus Z_j|+|W_j'|-n + |Z_i\setminus \tZ_i| - |W_{i}'|) \\
= \, &\sum_{i\in [r]}\sum_{j\in L\setminus \{i\}}|Z_j^i| (|\tZ_j\setminus Z_j|+|W_j'|-n + |Z_i\setminus \tZ_i|) - \sum_{i\in [r]} |\tZ_i\setminus Z_i| |W_{i}'| .
\end{aligned}
\]
Moreover, it is not hard to see that 
\[
\sum_{l=1}^{m} \sum_{j' <_L j} |Z_{j'}^i| = \sum_{i\in [r]}\sum_{\{j_1, j_2\}\in \binom{L\setminus \{i\}}2} |Z_{j_1}^i||Z_{j_2}^i| = \sum_{i\in [r]} e_T(\tZ_i\setminus Z_i).
\]
Now the claim follows by combining these estimates with \eqref{eq:eG0}.
\end{proof}

What remains is to estimate the number of edges in each $\tU_i$. 
For $i\in [r]$, we have $e(G[\tU_i]) = e(Z_i^i, G[\tU_i]) + e_G(\tU_i\setminus Z_i^i)$. 
To bound $e_G(\tU_i\setminus Z_i^i) = e_G((\tZ_i\setminus Z_i)\cup W_i')$, we note that $e_G(\tZ_i\setminus Z_i, W_i')\le |\tZ_i\setminus Z_i| |W_{i}'|$ and $e_G(W_i')\le e_T(W_i')$. 
However, we may not have $e_G(\tZ_i\setminus Z_i)\le e_T(\tZ_i\setminus Z_i)$ because each $Z_j^i$ is an independent set in $T$, but may not be independent in $G$ when $a\ge 2$.
Thus, $e_G(\tZ_i\setminus Z_i)\le e_T(\tZ_i\setminus Z_i)+\sum_{j\neq i}e_G(Z_j^i)$.
Putting these together, for each $i\in [r]$, we have
\begin{align}\label{eq:bat}
    e_G(\tU_i\setminus Z_i^i)=e_G(\tZ_i\setminus Z_i, W_i')+e_G(W_i')+e_G(\tZ_i\setminus Z_i) \le \alpha_i + \sum_{j\neq i}e_G(Z_j^i).
\end{align}
Let $f_i:= e(Z_i^i, G[\tU_i])$. 
By Claim~\ref{clm:eG'}, \eqref{eq:bat} and $e(G)= e(G') +  e(Z''\cup W'', G)  + \sum_{i\in [r]}e_G(\tU_i)$, we derive that
\begin{align}
e(G)&\le t_r(k) n^2 + e(Z''\cup W'', G) + \sum_{i\in [r]} \left(f_i + \beta_i -\alpha_i + e_G(\tU_i\setminus Z_i^i)\right) \label{eq:eG1} \\
&\le t_r(k) n^2 + e(Z''\cup W'', G) + \sum_{i\in [r]} \left(f_i + \beta_i + \sum_{j\ne i}e_G(Z_j^i)\right) \label{eq:eG2}
\end{align}

We now focus on the structure of each $\tU_i$.
We first show that $G[\tU_i]$ is ``almost'' $K_{t,t}$-free.

\begin{claim}
\label{clm:Ktt}
The following holds for all $i\in [r]$. 
\begin{enumerate}[label=$(K\arabic*)$]
\item Both $G[\tZ_i]$ and $G[Z_i^i\cup W_i']$ are $K_{t,t}$-free. \label{item:K1}
\item If $|W_i'|> 2t\eps n+2\gamma n$, then $|W_i'\setminus V_{q_i}|\le t-1$.
\label{item:K2}
\item If $|W_i'|> 2t\eps n+2\gamma n$, then $G[\tZ_i\cup (W_i'\cap V_{q_i})]$ is $K_{t,t}$-free.\label{item:K3}
\end{enumerate}
\end{claim}

\begin{proof}
For~\ref{item:K1}, suppose there is a copy of $K_{t,t}$ in $\tU_i$, with vertex set denoted by $B$, contained in $\tZ_i$ or in $Z_i^i\cup W_i'$.
Let $N_B$ be the set of common neighbors of these $2t$ vertices of $B$.
First assume that $B\subseteq \tZ_i$.
Then for any $j\in L\setminus\{i\}$, by~\ref{item:DegZ} we have $|N_B\cap W_j'|\ge |W_j'| - 4t\eps n$, and thus by~\ref{item:Wsizes} $|N_B\cap W_j\cap W_j'|\ge |W_j'| - 4t\eps n - 2\gamma n \ge n/2$.
For any $j\notin L\cup\{i\}$, because $B\cap Z_j=\emptyset$ by~\ref{item:WiLarge}, we have $|N_B\cap Z_j|\ge an - 4t\eps n\ge n/2$ by~\ref{item:DegZ}.
Note that every set in $\{N_B\cap Z_j: j\notin L\}\cup \{N_B\cap W_j\cap W_j': j\in L\}$ has size at least $n/2$ and every pair of them has density at least $1-4\gamma$.
%
Therefore we can find a copy of $K_{r-1}{(t)}$ in the union of these sets, which gives rise to a copy of $K_{r+1}(t)$ together with $B$, a contradiction.

Second we assume that $B\subseteq Z_i^i\cup W_i'$.
In this case we note that for any $j\neq i$, we have $B\cap Z_j=\emptyset$ and thus by~\ref{item:DegZ}, we have $|N_B\cap Z_j^j|\ge (1-\sqrt\gamma)an - 4t\eps n\ge n/2$.
Then as these sets have high pairwise densities, as in the previous case, we can find a copy of $K_{r-1}{(t)}$ in the union of these sets, yielding a copy of $K_{r+1}(t)$ together with $B$, a contradiction.
Now~\ref{item:K1} is proved.

Now we turn to~\ref{item:K2}, and suppose $|W_i'|> 2t\eps n+2\gamma n$ and thus $|W_i\cap W_i'|> 2t\eps n$ by~\ref{item:Wsizes}.
First, if $W_i'$ contains at least $t$ vertices which are not from $V_{q_i}$ (the cluster containing $W_i$), then by~\ref{item:DegZ}, each of these vertices have at most $2\eps n$ non-neighbors in $W_i\cap W_i'$, and thus we can find a copy of $K_{t,t}$ in $W_i'$, contradicting~\ref{item:K1}.
So we have $|W_i'\setminus V_{q_i}|\le t-1$.

For~\ref{item:K3}, suppose there is a copy of $K_{t,t}$ as stated in the claim, whose vertex set is denoted by $B$.
As in the previous paragraph, we have $|W_i|> 2t\eps n$ by~\ref{item:Wsizes}.
Now observe crucially that if $B\cap Z_j^i\neq \emptyset$, then by~\ref{item:WiLarge} $|W_i|+|W_j| > n$, and thus, $W_i$ and $W_j$ are not from the same cluster.
So by~\ref{item:DegZ}, for any $j\in [r-1]\setminus \{i\}$, if $B\cap Z_{j}^i=\emptyset$, then the vertices of $B$ have large common neighborhoods in $Z_{j}^j$; if $B\cap Z_{j}^i\neq \emptyset$, then the vertices of $B$ have large common neighborhoods in $W_j\cap W_j'$ (note that $|W_j|\ge (1-\eps-\sqrt\gamma)n$ by~\ref{item:WiLarge}).
Since each of these common neighborhoods have size at least $n/2$ and each pair of them has high density, we can find a copy of $K_{r-1}{(t)}$ in the union of these sets, yielding a copy of $K_{r+1}(t)$ together with $B$, a contradiction.
%
%
\end{proof}

%

We now derive a lower bound for $\sum f_i$ from Claims~\ref{clm:eG'} and \ref{clm:Ktt}. For $i\in [r]$, we have $\beta_i \le \sum_{j\in L\setminus \{i\}}|Z_j^i| (|\tZ_j\setminus Z_j|+ |W_j'\setminus V_{q_j}| + |Z_i\setminus \tZ_i|)$ as $|W_j'|-n\le |W_j'\setminus V_{q_j}|$.
Fix $j\in L\setminus \{ i \}$. Note that $|W_j|\ge (1-2\eps)n$.
We have $|\tZ_j\setminus Z_j|\le t-1$ by~\ref{item:Wismall}, and $|W_j'\setminus V_{q_j}|\le t-1$ by~\ref{item:K2}.
If $|W_i| > n/2$, then $|Z_j^i|\le t-1$ by~\ref{item:Wismall}. 
Furthermore, since $|Z_i\setminus \tZ_i|\le (r-1)\sqrt\gamma n + C_0$ by \ref{item:Zsizes}, it follows that
\[
|Z_j^i| \left(|\tZ_j\setminus Z_j|+|W_j'\setminus V_{q_j}| + |Z_i\setminus \tZ_i|\right)
\le (t-1) (t-1 + t-1 + (r-1)\sqrt\gamma n + C_0)\le (t-1)r\sqrt\gamma n.
\]
Otherwise $|W_i| \le n/2$, and by~\ref{item:WiLarge}, we have $Z_i^{i'}=\emptyset$ for any $i'\neq i$.
This implies $|Z_i\setminus \tZ_i|=|Z_i''|\le C_0$. Using $|Z_j^i|\le \sqrt\gamma n$, \ref{item:Wismall}, and \ref{item:K2}, we derive that  
\[
|Z_j^i| \left(|\tZ_j\setminus Z_j|+|W_j'\setminus V_{q_j}| + |Z_i\setminus \tZ_i|\right)\le \sqrt\gamma n (2(t-1)+ C_0)\le 2C_0\sqrt\gamma n.
\]
Summarizing these two cases for all $j\in L\setminus \{i\}$, we obtain that $\beta_i\le (r-1) 2C_0\sqrt\gamma n$, and consequently,
\begin{equation}
\label{eq:bi}
\sum_{i\in [r]}\beta_i \le 2(r-1) r C_0\sqrt\gamma n.
\end{equation}

On the other hand, for all $i\ne j$, the graph $G[Z_j^i]$ is $K_{t, t}$-free by~\ref{item:K1} and thus, by \ref{item:Zsizes}, $\sum_{i, j: i\neq j}e_G(Z_j^i)\le r(r-1) C \left( \sqrt\gamma n\right)^{2-1/t}$.
Applying this with \eqref{eq:eG2}, \eqref{eq:bi}, and the fact that $e(Z''\cup W'', G)\le (r+1)C_0 kn$, we obtain that
\begin{align*}
    e(G)&\le t_r(k) n^2 + (r+1)C_0 kn + \sum_{i\in [r]}f_i +  2(r-1) r C_0\sqrt\gamma n + \sum_{i, j: i\neq j}e_G(Z_j^i)  \\
    &\le t_r(k) n^2 + \sum_{i\in [r]}f_i + r^2 C\sqrt{\gamma} n^{2-1/t},
\end{align*}
as $\gamma \ll 1$.
Using the assumption $e(G)\ge g(n,r,k,t)\ge t_r(k) n^2 + z_{t}(n)$, we infer that
\begin{align}
\label{eq:edge-rough}
\sum_{i\in [r]}f_i \ge z_{t}(n) - r^2 C\sqrt{\gamma} n^{2-1/t} \ge \frac{c}2 n^{2 - 1/t}
\end{align}
by using \ref{item:zest},  $z_t(n)\ge z_t^{(2)}(n) \ge c n^{2 - 1/t}$, and $\gamma \ll 1$.

\medskip
We next study the existence of $K_{1,t}$ in each color class. To do so,
we consider a copy of $K_3(t)$ in $G[\tU_i\cup \tU_j]$ for some $i\ne j$.
\begin{claim}
\label{clm:K3t}
For any $i\ne j$, if $G[\tU_i\cup \tU_j]$ contains a copy $K$ of $K_3{(t)}$, then there exists $l\notin \{i, j\}$ such that $V(K)$ intersects $V_{q_{l}}$ and every cluster in $Z_{l}$.
\end{claim}

\begin{proof}
We may assume that $r>2$ as otherwise the claim is trivial.
Suppose to the contrary that there is a copy $K$ of $K_3{(t)}$ in, say, $\tU_1$ and $\tU_2$, such that for every $l\in [3, r]$, there is a cluster in $U_l$ which does not intersect $B:=V(K)$.
Let $V_{i_l}$ be a cluster in $Z_l$ such that $B\cap V_{i_l}=\emptyset$, and if there is no such cluster in $Z_l$, then we choose $V_{i_l}=V_{q_l}$.
Note that in the former case, we have $|\tU_l\cap V_{i_l}|=|Z_l^l\cap V_{i_l}|\ge (1-\sqrt\gamma a)n$.
In the latter case, we have $Z_l^1\neq \emptyset$ or $Z_l^2\neq \emptyset$, which implies that $|W_{l}|\ge (1-2\eps)n$ by~\ref{item:WiLarge}, and thus $|\tU_l\cap V_{i_l}|=|W_l'\cap V_{i_l}|\ge (1-3\eps)n$.
Now, by~\ref{item:DegZ}, every vertex in $B$ has at most $2\eps n$ non-neighbors in $\tU_l\cap V_{i_l}$ for each $l\in [3, r]$.
Since for every $l$ we have $|\tU_l\cap V_{i_l}|\ge 0.9n$, 
%
one can find large common neighborhoods (e.g. of size $n/2$) of all vertices of $B$ in each $\tU_l\cap V_{i_l}$, and then find a copy of $K_{r-2}{(t)}$ in these sets.
Altogether we obtain a copy of $K_{r+1}(t)$, a contradiction.

Therefore, for such a copy $K$ of $K_3(t)$, there exists $l\notin \{i,j\}$ such that $K$ must intersect all clusters of $U_l$.
Since $V(K)\cap Z_l\neq \emptyset$, we have $Z_l^i\neq \emptyset$ or $Z_l^j\neq\emptyset$.
Then by~\ref{item:WiLarge}, $|W_{l}|\ge (1-2\eps)n$ and in particular, $V_{q_l}\neq\emptyset$. Therefore $V(K)\cap V_{q_l}\neq \emptyset$.
\end{proof}

\begin{claim}
\label{clm:K1t}
For all but exactly one $j\in [r]$, we have $d(v, Z_j^j)\le t-1$ for all $v\in \tU_j$.
\end{claim}

\begin{proof}
First assume that there exists $j\in [r]$ such that $G[\tU_j]$ contains a copy of $K_{1,t}$, with vertex set denoted by $\{v, u_1,\dots, u_t\}$, $v\in \tU_j$ and $u_1,\dots, u_t\in Z_j^j$.
Fix $i\in [r]\setminus \{j\}$ and let $N'$ be the set of common neighbors of $u_1,\dots, u_t$ in $\tU_i\cap U_i$.
Suppose $v\in V_p$ and let $N$ be the set of common neighbors of these $t+1$ vertices in $\tU_i\cap U_i$.
In particular, $N\subseteq N'$ and $N$ is almost equal to the union of $a$ or $a+1$ clusters in $\tU_i$.
Suppose there is a copy of $K_{t-1,t}$ with parts $S_1$ of size $t-1$ and $S_2$ of size $t$ such that $S_1\subseteq N'$ and $S_2\subseteq N$.
Then by Claim~\ref{clm:K3t}, there exists $l\in [r]\setminus \{i,j\}$ such that $B\cap Z_l\neq\emptyset$ and $B\cap V_{q_l}\neq\emptyset$, where $B$ denotes the vertex set of the copy of $K_3{(t)}$.
This is impossible since $v$ is the only possible vertex in $B\cap (Z_l\cup V_{q_l})$ and can not satisfy both.
Therefore, 
letting $N^*=N\cup (N'\cap V_p)$, we infer that $e_G(N^*)=e_G(N)+e_G(N, N'\setminus N)=O(n^{2-1/(t-1)})$.

By~\ref{item:Wsizes},~\ref{item:Zsizes} and~\ref{item:DegZ}, we have $|\tU_i\setminus N^*|\le 3(t+1)\eps n$.
Let $E^i$ be the set of the edges incident to $\tU_i\setminus N^*$ and counted in $f_i$.
We split it to 
$E^i\cap E_G(Z_i^i)$ and $E^i\cap E_G(\tU_i\setminus Z_i^i, Z_i^i)$.
Note that by~\ref{item:K1}, each of the terms can be split further into at most $k$ $K_{t,t}$-free bipartite graphs, each with one part of size at most $3(t+1)\eps n$ and the other part of size at most $(1+(r-2)\sqrt\gamma)an$.
Therefore, we obtain that
\begin{equation}
\label{eq:U-W}
f_i= O(\eps n^{2-1/t}) + O(n^{2-1/(t-1)}) =O(\eps n^{2-1/t}).
\end{equation}

Now assume there exist distinct $j_1, j_2\in [r]$ such that each $G[\tU_{j_i}]$ contains a copy of $K_{1,t}$ whose part of size $t$ is in $Z_{j_i}^{j_i}$.
The arguments above imply that \eqref{eq:U-W} holds for all $i\in [r]$, and consequently, $\sum_{i\in [r]}f_i= O(\eps n^{2-1/t})$, contradicting~\eqref{eq:edge-rough}.  

On the other hand, if $d(v, Z_j^j)\le t-1$ for all $j\in [r]$ and all $v\in \tU_j$, 
then $\sum_{j\in [r]} f_j\le (t-1) kn$,  
again contradicting \eqref{eq:edge-rough}. 
\end{proof}

By Claim~\ref{clm:K1t}, without loss of generality, we assume that, 
\begin{align}\label{eq:toZii}
  \text{for } i\ge 2, \quad d(v, Z_i^i)\le t-1 \ \text{for all} \ v\in \tU_i, \quad \text{and thus}, \quad f_i\le \begin{cases}
      (t-1)|\tU_i| \hfill\text{ if } a\ge 2,\\
      (t-1)|\tU_i\setminus Z_i^i| \hfill \text{ if } a=1.
  \end{cases}
\end{align}


If $b=0$, then $\tU_i = Z_i^i = Z_i \setminus Z''$ for all $i$ by \eqref{eq:Ui0}. 
In this case $\tU_1$ is $K_{t,t}$-free by~\ref{item:K1} and $\tU_i$ is $K_{1,t}$-free for all $i\ge 2$ by \eqref{eq:toZii}. 
Since $G$ is $\gamma$-close to $K_r(an)$, $G[U_i\setminus Z'', U_j\setminus Z'']$ is almost complete for all $i\ne j$. 
This completes the proof of Theorem~\ref{thm7} with $Z:=Z''$. \qed


\bigskip
By~\eqref{eq:edge-rough} and \eqref{eq:toZii}, we get
\begin{equation}
\label{eq:f1x}
 f_1\ge z_{t}(n) - \eps n^{2-1/t}.
 \end{equation}
In particular, we claim that
\begin{equation}
\label{eq:W1x}
|W_1|> 3t\eps n \quad \text{ if } b>0
\end{equation}
(which we will refine a moment later).
Indeed, the edges counted in $f_1$ can be covered by $G[Z_1^1]$, $G[Z_1^1, W_1\cap W_1']$, and at most $k$ $K_{t,t}$-free bipartite graphs, each with a part of size at most $\sqrt\gamma n$ and a part of size at most $an$.
If $|W_1|\le 3t\eps n$, then $e_G(Z_1^1, W_1\cap W_1') = O(\eps n^{2-1/t})$.
Together with $e_G(Z_1^1)\le z_{t}^{(a)}(n)$,
we have
\[
f_1\le z_{t}^{(a)}(n) + O(\eps n^{2-1/t}) < z_{t}^{(a+1)}(n) - \eps n^{2-1/t}
\]
by~\ref{item:za=1gap}, contradicting~\eqref{eq:f1x}.

\medskip
Now we can give a much cleaner structure, shown in a series of claims below.
A key step is to show that $Z''\cup W''=\emptyset$.
{From now on we only consider $k\le 2r$}.

\begin{claim}\label{clm:ZW''}
Suppose $v_0\in V(G)$ and $i\in [r]$ satisfy that $v_0$ has at least $\eps n$ neighbors in $Z_j$ for every $j\neq i$.
Then $v_0$ has less than $\eps n$ neighbors in $U_i$.
In particular, we have $Z''=\emptyset$ and $W''=\emptyset$. 
\end{claim}

\begin{proof}
The second part of the claim follows immediately from the definitions of $Z''$ and $W''$.

Suppose to the contrary, that there exist $v_0\in V(G)$ and $i\in [r]$ such that $v_0$ has at least $\eps n$ neighbors in $Z_j$ for every $j\neq i$ and at least $\eps n$ neighbors in $U_i$.
Since $|Z_j^j | \ge (1-\sqrt\gamma)an$ for all $j \in [r]$, there exist sets $N_1,\dots, N_{r-1}$ each of size $\eps n-\sqrt\gamma n$ such that $N_j\subseteq Z_j^j\cap N(v_0)$ for $j\neq i$ and $N_i\subseteq (Z_i^i\cup W_i)\cap N(v_0)$.
Recall that $W_1'=W_1'\cap V_{q_1}$.
By averaging, there exists $N_1'\subseteq N_1$ with $|N_1'|\ge (\eps n -\sqrt\gamma n - 2r\gamma n)/2 \ge \eps n/3$ such that all vertices of $N_1'$ are in $Z_1^1\cup W_1'$ and from the same cluster, that is, 
\[
N_1'\subseteq Q, \text{ where } Q\in \{ Z_1^1, W_1'\}.
\]
Note that $N_1'\subseteq W_1'$ is possible only if $i=1$ and $a=1$.
If $i\neq 1$, then let $N_i':=N_i\setminus ((W_i\setminus W_i')\cup V_{q_1})$ and for every $j\in [r]\setminus \{1, i\}$, let $N_j':=N_j$.
By~\ref{item:Wsizes}, $|W_i\setminus W_i'|\le 2r\gamma n$, and by~\eqref{eq:W1}, $|W_i\cap V_{q_1}|\le \gamma n$.
Thus, we have $|N_j'|\ge \eps n/3$ for all $j\in [r]$.
Because the sets $N_j'$ are small, we can not apply the degree conditions~\ref{item:DegZ} to them and instead, we use~\ref{item:Deg0}.

Recall that $B_1$ is given by~\ref{item:Deg0}.
Next we show that $G[\tU_1\setminus B_1]$ does not contain a copy of $K_{t-1, t}$ such that the part of size $t$ is in $N_1'$.
Suppose instead, there is such a copy of $K_{t-1, t}$, with parts denoted by $A$ and $B$, such that $|A|=t$, $A\subseteq N_1'\setminus B_1$ and $B\subseteq \tU_1\setminus B_1$. 
Recall that $N_i'\cap V_{q_1}=\emptyset$ and for each $j\in [r]\setminus \{1,i\}$, $N_j'\subseteq Z_j^j$.
Observe that for every $v\in \tU_1\setminus B_1$, we have $d(v, N_j')\ge |N_j'| - \sqrt\gamma n$.
Indeed, if $j\neq i$, then $N_j'\subseteq Z_j^j$ and we have $d(v, N_j')\ge |N_j'| - \sqrt\gamma n$ by~\ref{item:Deg0}; otherwise note that $N_i' \subseteq Z_i^i\cup (W_i'\cap W_i)$, and by~\ref{item:Deg0} and $N_i'\cap V_{q_1}=\emptyset$ we have $d(v, N_i')\ge |N_i'| - \sqrt\gamma n$.
Therefore, we obtain that the vertices in $A\cup B$ have at least $|N_j'| - (2t-1)\sqrt\gamma n\ge (1- \gamma^{1/3})|N_j'|$ common neighbors in each $N_j'$, $j\in [2,r]$.
Because each pair $N_j', N_{j'}'$ has a high density, we can find a copy of $K_{r-1}{(t)}$ in the union of these common neighborhoods, which together with $A\cup B\cup \{v_0\}$ form a copy of $K_{r+1}(t)$, a contradiction.

Now given that $G[\tU_1\setminus B_1]$ does not contain a copy of $K_{t-1, t}$ such that the part of size $t$ is in $N_1'\setminus B_1$, we give a refined estimate on $f_1$. 
Indeed, since $G[N_1'\setminus B_1, Z_1^1\setminus B_1]$ does not contain a copy of $K_{t-1, t}$ such that the part of size $t$ is in $N_1'\setminus B_1$, we get $e_G(N_1'\setminus B_1, Z_1^1\setminus B_1)= O(n^{2-1/(t-1)})$.
Similarly $e_G(N_1'\setminus B_1, W_1'\setminus B_1)= O(n^{2-1/(t-1)})$.
Suppose $N_1'\subseteq V_q$ for some $q\in \{1, q_1\}$, 
then we have
\[
E(G[\tU_1]) = E(G[\tU_1\setminus (N_1'\setminus B_1)]) \cup E(G[N_1'\setminus B_1, \tU_1\setminus (B_1\cup V_q)]) \cup E(G[N_1'\setminus B_1, B_1\cap \tU_1)]).
\]
Recall that $|N_1'|\ge \eps n/3$ and $|B_1|\le 2\sqrt\gamma n$.
Therefore, (regardless of $a=1$ or $(a,b)=(2, 0)$) we can bound $f_1\le |E(G[\tU_1])|$ by
\[
f_1 \le z_t\left((1-\tfrac{\eps}{3})n,n\right) + O(n^{2-1/(t-1)})+ O(\sqrt\gamma n^{2-1/t}) <  z_t(n) -  3rC_0kn,
\]
where we used~\ref{item:zagap} and $\gamma \ll \eps$.
This contradicts~\eqref{eq:f1x}.
\end{proof}

{When $a=2$ and $b=0$ (i.e., $k=2r$), since $Z''_i=\emptyset$ and $W_i=\emptyset$ for all $i\in [r]$, 
by~\eqref{eq:Ui0}, we get $\tU_i=Z_i$ for all $i\in [r]$.
Therefore 
$e(G) = e(G') + \sum_{i=1}^r e_G(\tU_i) + e(Z''\cup W'', G)\le t_r(k)n^2 + z_t(n) + (r-1)(t-1)n$ by~\ref{item:K1} and~\eqref{eq:toZii}, proving Theorem~\ref{thm:lb} for $k=2r$.
}

\bigskip
For the remaining of the proof, we only need to consider $a=1$ (and thus $b>0$).
Moreover, now for $i,j\in [r]$ each $Z_i^j\subseteq Z_i$ is an independent set and thus $e_G(Z_i^j)=0$.
So we can first update our bounds on $e(G)$ and $f_1$.
Recall the bounds~\eqref{eq:eG2},~\eqref{eq:bi} and~\eqref{eq:toZii} and we have
\begin{align*}
    e(G)&\le t_r(k) n^2 + \sum_{i\in [r]}(f_i+\beta_i)    \\
    &\le t_r(k) n^2 + f_1 + (r-1)(t-1)(1+\sqrt{\gamma})n + 2(r-1) r C_0\sqrt\gamma n,
\end{align*}
yielding
\begin{equation}
\label{eq:f1refined}
    f_1 \ge z_t(n,n) - C_0 n
\end{equation}

\begin{claim} 
\label{clm:clean}
Suppose $b>0$. Then $\tU_1=Z_1^1\cup W_1'$ and $W_1'\subseteq V_{q_1}$.
\end{claim}

\begin{proof}
Suppose to the contrary, there is a vertex $v$ in $\tU_1\setminus (Z_1^1\cup W_1')$ or $W_1'\setminus V_{q_1}$, namely, $v\in Z_i^1$ for some $2\le i\le r$ or $v\in W_1'\setminus V_{q_1}$.
Suppose $v\in V_{l}$. Then $l\ne q_1$.
Moreover, if $i$ is defined, then $W'_1\cap V_{q_1}\subseteq V\setminus (\tZ_i \cup V_{l})$; otherwise, $W'_1\cap V_{q_1}\subseteq V\setminus V_{l}$.
By 
\ref{item:DegZ}, we have $\overline d(v, W_1'\cap V_{q_1}) \le 2\eps n$.
Let $N:=W_1'\cap V_{q_1}\cap N(v)$. We have $|(W_1'\cap V_{q_1})\setminus N| \le 2\eps N$.
Since $|W_1'\setminus V_{q_1}| \le |W_1'\setminus W_1|\le 2\gamma n$, it follows that $|W'_1\setminus N|\le 2\eps n + 2\gamma n \le 3\eps n$.

Recall \eqref{eq:W1x}, $|W_1|> 3t\eps n$. By~\ref{item:K3} (if $v\in \tZ_1\setminus Z_1$) or~\ref{item:K1} (if $v\in W_1'\setminus V_{q_1}$), we know that $G[Z_1^1, N]$ contains no $K_{t-1, t}$ with the part of size $t$ in $N$.
This implies that $e_G(Z_1^1, N)= O(n^{2-1/(t-1)})$.
Furthermore, by \ref{item:Zsizes} and \ref{item:K1}, $G[\tZ_1\setminus Z_1^1, Z_1^1]$ is a $K_{t,t}$-free bipartite graph with one part of size at most $(r-1)\sqrt{\gamma} n$ and the other part of size at most $n$. Thus, $e_G(\tZ_1\setminus Z_1^1, Z_1^1) \le C(r-1)\sqrt{\gamma} n^{2-1/t}$. By the similar arguments, we have $e_G( W_1'\setminus N, Z_1^1) \le C (3\eps n) n^{1-1/t}$.

Putting these bounds together (and note that $Z_1^1$ is an independent set, we get
\begin{align*}
    f_1 &= e_G(Z_1^1, N) + e_G(\tZ_1\setminus Z_1^1, Z_1^1) + e_G( W_1'\setminus N, Z_1^1)\\
    & = O(n^{2-1/(t-1)}) + O(\sqrt{\gamma} n^{2-1/t}) + O(\eps n^{2-1/t}).
\end{align*}
By~\ref{item:za=1gap}, this contradicts~\eqref{eq:f1refined}.
\end{proof}

Claim~\ref{clm:clean} shows that $\tU_1$ has no atypical vertices and is thus $K_{t,t}$-free by~\ref{item:K1}. 
Furthermore, since $\tU_1=Z_1^1\cup W_1'$ and $W_1'\subseteq V_{q_1}$, it follows that
\begin{align}
    \label{eq:etU1}
    \alpha_1 = \beta_1 = 0, \quad \text{and} \quad
    e_G(\tU_1) = f_1 \le z_t (|Z_1^1|, |W'_1|). 
\end{align}
Therefore, if $|W_1'|\le (1-\gamma)n$, then we have $f_1\le z_{t}(n, |W_1'|) \le z_{t}(n,n) - \delta n^{2-1/t}$ for some $\delta>0$
by~\ref{item:zagap}.
This contradicts~\eqref{eq:f1refined}.
So we obtain
\begin{equation}
\label{eq:W1}
\text{if } a=1, \text{ then }|W_1'|\ge (1-\gamma)n \quad  (\text{and thus} \quad 1\in L).
\end{equation}

Next we study $G[\tU_i]$ for $i\ge 2$.
A key observation is that copies of $K_{1,t}$ in $G[\tU_i]$ together with copies of $K_{t-1,t}$ in $\tU_1$ may form copies of $K_3{(t)}$, which are restricted by Claim~\ref{clm:K3t}.

\begin{claim} 
\label{clm:clean0}
Suppose $i\in [2, r]$. 
\begin{enumerate}
\item If there is a copy of $K_{1,t}$ in $\tU_i\setminus (Z_1\cup V_{q_1})$, then there exists $l\in [r]\setminus \{i\}$ such that the vertex set of $K_{1,t}$ intersects both $V_{q_l}$ and $Z_l$.
\item Both $\tZ_i\setminus Z_1$ and $Z_i^i \cup (W_i'\setminus V_{q_1})$ are $K_{1,t}$-free.
\end{enumerate}
\end{claim}

\begin{proof}
For Part (1), let $B$ be the vertex set of a copy of $K_{1,t}$ in $\tU_i\setminus (Z_1\cup V_{q_1})$. Since $B\cap (Z_1\cup V_{q_1})=\emptyset$ and $\tU_1\subseteq Z_1\cup V_{q_1}$, by~\ref{item:DegZ}, all vertices of $B$ have at most $2\eps n$ non-neighbors in $\tU_1$.
Letting $N:=\tU_1\cap \bigcap_{w\in B}N(w)$, we have $|N|\ge |\tU_1| - (t+1)2\eps n$.

First assume that $N$ is $K_{t-1,t}$-free and thus $e_G(N) = O(n^{2-1/(t-1)})$.
Note that, since $|\tU_1\setminus N|\le (t+1)2\eps n$, the edges in $\tU_1$ incident to $\tU_1\setminus N$ can be split into two bipartite $K_{t,t}$-free graphs each with one part of size at most $(t+1)2\eps n$ and the other part of size at most $n$. 
Thus, the number of such edges is $O(\eps n^{2-1/t})$.
This gives $f_1 =  O(n^{2-1/(t-1)})+O(\eps n^{2-1/t})$, contradicting~\eqref{eq:f1x}.

We thus assume $N$ contains a copy of $K_{t-1,t}$. Together with $B$, they form a copy of $K_3{(t)}$ in $G[\tU_1\cup \tU_i]$ and we denote its vertex set by $B'$.
By Claim~\ref{clm:K3t}, there exists $l\notin\{1,i\}$ such that $B'$ intersects $V_{q_l}$ and $Z_l$.
By Claim~\ref{clm:clean}, $\tU_1\cap U_l=\emptyset$, so $B'\cap Z_l=B\cap Z_l$ and $B$ indeed intersects $Z_l$.
Since $\tU_i\cap Z_l\supseteq B\cap Z_l\neq \emptyset$, we infer that $|W_l|\ge (1-2\eps)n$ from \ref{item:Wismall}, which implies that $q_l\neq q_1$ because of \eqref{eq:W1x}.
It follows that $W_1\cap V_{q_l}=\emptyset$ and thus $B\cap V_{q_l} = B'\cap V_{q_l} \neq \emptyset$, as desired.

For Part (2), let $A_i:=Z_i^i \cup (W_i'\setminus V_{q_1})$ and 
$B$ be the vertex set of a copy of $K_{1,t}$ in $\tZ_i\setminus Z_1$ or in $A_i$.
Then, by the first part of the claim, there exists $l\in [r]\setminus \{i\}$ such that $B$ intersects $V_{q_l}$ and $Z_l$.
This is impossible if $B\subseteq A_i$ because $A_i\cap Z_l=\emptyset$ for any $l\notin \{1, i\}$, and also impossible if $B\subseteq \tZ_i\setminus Z_1$ because in which case $B\cap W=\emptyset$ and thus $B\cap V_{q_l}= \emptyset$ for any $l\notin \{1, i\}$.
%
\end{proof}

The following claim shows a clean structure for the $\tU_i$ such that $W_i'$ is  not too small.

\begin{claim}
\label{eq:Wiepsn}
For $i\in [2,r]$ such that $|W_i|\ge 2\eps n$, we have $\tU_i\subseteq U_i\cup V_{q_i}$.
\end{claim}

\begin{proof}
Suppose instead, for some $i_0\in [2,r]$ with $|W_{i_0}|\ge 2\eps n$, there exists $v\in \tU_{i_0}\setminus (U_{i_0}\cup V_{q_{i_0}})$.
{By \ref{item:Zpartition} and the fact that $v\in \tU_{i_0}\setminus U_{i_0}$, we infer that $d(v, Z_j)\ge \eps n$ for all $j\neq {i_0}$.
Then, by Claim~\ref{clm:ZW''}, we have $d(v, U_{i_0})<\eps n$.
Consequently, $d(v, W_{i_0}'\cap W_{i_0}) < \eps n$, namely, $v$ has at least $2\eps n - 2\gamma n - \eps n\ge (1/2)\eps n$ non-neighbors in $W_{i_0}'\cap W_{i_0}$ (in $G$).}
Note that $v$ is adjacent to all the vertices of $W_{i_0}'\cap W_{i_0}$ in $T$.
Since $G[W'_{i_0}]\subseteq T[W'_{i_0}]$, we infer that 
\[
e_G(W_{i_0}') + e_G(\tZ_{i_0}\setminus Z_{i_0}, W_{i_0}') \le e_T(W_{i_0}') + |\tZ_{i_0}\setminus Z_{i_0}||W_{i_0}'| - (1/2)\eps n.
\]
Since $e_G(\tZ_{i_0}\setminus Z_{i_0})\le e_T(\tZ_{i_0}\setminus Z_{i_0})$
and $\alpha_{i_0}=|\tZ_{i_0}\setminus Z_{i_0}| |W_{i}'| +e_T(W_{i_0}')+e_T(\tZ_{i_0}\setminus Z_{i_0})$, we have
\begin{align}
e_G(\tU_{i_0}\setminus Z_{i_0}) &= e_G(\tZ_{i_0}\setminus Z_{i_0}) + e_G(W_{i_0}') + e_G(\tZ_{i_0}\setminus Z_{i_0}, W_{i_0}') \nonumber\\ 
&\le e_T(\tZ_{i_0}\setminus Z_{i_0}) + e_T(W_{i_0}') + |\tZ_{i_0}\setminus Z_{i_0}||W_{i_0}'| - (1/2)\eps n \nonumber\\
&\le \alpha_{i_0} - (1/2)\eps n. \label{eq:i0}
\end{align}
Combining \eqref{eq:bat} and \eqref{eq:i0} gives  
\begin{align}\label{eq:eUZ}
    \sum_{i\in [r]} e_G(\tU_i\setminus Z_i) \le \sum_{i\in [r]} \alpha_i - (1/2)\eps n
\end{align}
 
Recall that $f_i\le (t-1)|\tU_i\setminus Z_i|$ ($i\ge 2$) by~\eqref{eq:toZii}, $|W'_1|\ge (1-\gamma)n$ by~\eqref{eq:W1}, and $|Z_i^i|\ge (1-\sqrt\gamma)n$ by~\ref{item:Zsizes}.
Therefore, as $\sum_{i=1}^r |\tU_i\setminus Z_i|\le (b-1)n + r\sqrt{\gamma}n +   {\gamma}n$, we obtain
\begin{align}\label{eq:fi}
 \sum_{i=2}^r f_i\le (t-1) \left((b-1)n + r\sqrt{\gamma}n +   {\gamma}n \right) \le (t-1)(b-1)n + \sqrt[3]{\gamma}n.
\end{align}

Since $Z''\cup W''=\emptyset$, \eqref{eq:eG2} becomes
$e(G) \le t_r(k)n^2 +\sum_{i\in [r]}(f_i+\beta_i  + e_G(\tU_i\setminus Z_i)- \alpha_i )$. 
Recall that $\sum_{i=1}^r \beta_i\le 2r^2 C_0\sqrt\gamma n$ by \eqref{eq:bi}.
Together with \eqref{eq:eUZ} and \eqref{eq:fi}, we derive that
\[
e(G)\le t_r(k)n^2 + z_t( n) + (t-1)(b-1)n + \sqrt[3]\gamma n + 2r^2 C_0\sqrt\gamma n -\eps n/2 < g(n,r,k,t),
\] 
as $\gamma \ll \eps$.
This is a contradiction.
\end{proof}

Let $L_1\cup L_2\cup L_3$ be a partition of $[2, r]$ such that $i\in L_1$ if and only if $|\tZ_i|<n$, $i\in L_2$ if and only if $|\tZ_i|=n$, and $i\in L_3$ if and only if $|\tZ_i|>n$.
The following properties hold for $L_1$, $L_2$ and $L_3$.
\begin{enumerate}[label=$(R\arabic*)$]
\item If $i\in L_1$, then $Z_i^j \ne \emptyset$ for some $j\ne i$. By \ref{item:WiLarge}, we have $i\in L$ and, by Claim~\ref{eq:Wiepsn}, $\tZ_i=Z_i^i\subsetneq V_i$ and $W_i'\subseteq V_{q_i}$.  \label{item:R1}
\item If $i\in L_2$, then $\tZ_i=Z_i^i=V_i$. Indeed, otherwise $\tZ_i\ne Z_i^i$, then $|Z_i^i|< n$ and $Z_i^j \ne \emptyset$ for some $j\ne i$. 
By~\ref{item:WiLarge} and Claim~\ref{eq:Wiepsn}, we have $\tZ_i=Z_i^i$, a contradiction. 
\item If $i\in L_3$, then $\tZ_i \not\subseteq Z_i$ (otherwise $|\tZ_i|\le n$). By Claim~\ref{eq:Wiepsn}, we have $|W_i|< 2\eps n$, which implies that $Z_i^j= \emptyset$ for $j\ne i$ by \ref{item:WiLarge}. Thus, $Z_i^i= Z_i= V_i\subsetneq \tZ_i$.
\end{enumerate}


Now we derive our final bound on $e(G)$. Write $z_i:=|\tZ_i|$ and $w_i:=|W_i'|$ for $i\in [r]$. 

By Claim~\ref{clm:eG'} and the fact that $Z''\cup W''=\emptyset$, we have
\[
e(G)=e(G')+ \sum_{i\in [r]}e_G(\tU_i) \le t_r(k) n^2 + \sum_{i\in [r]}\left(\beta_i-\alpha_i + e_G(\tU_i)\right).
\]
Moreover, as $a=1$, \eqref{eq:bat} becomes $e_G(\tU_i\setminus Z_i^i)\le \alpha_i$. It follows that $e_G(\tU_i) = f_i +   e_G(\tU_i\setminus Z_i) \le f_i + \alpha_i$.
For $i\in \{1\}\cup L_1\cup L_2$, we simply use $f_i$ as the upper bound and thus we get
\[
e_G(\tU_i)-\alpha_i\le f_i\le \begin{cases}
z_t(z_1, w_1) \quad \hfill \text{ if } i=1 \text{ by ~\eqref{eq:etU1},}\\
(t-1)\min\{z_i, w_i\} \quad \hfill \text{ if } i\in L_1 \text{ by $(R1)$ and Claim~\ref{clm:clean0} (2),} \\
(t-1) (z_i-n+w_i) \quad \hfill \text{ if } i\in L_2 \text{ by $(R2)$ and~\eqref{eq:toZii}.}
\end{cases}
\]
Additional work is needed for $i\in L_3$.
We let $\lambda = \max\{0, |L_1|-|L_3|\}$ and for $i\in L_3$, let $\lambda_i$ be the number of indices $j\in L_1$ such that $Z_j^i\neq \emptyset$.
By $(R1)$--$(R3)$, we know that if $Z_j^i\ne \emptyset$, then $i\in L_3$ and $j\in L_1$. 
This implies that $\lambda_i\ge 1$ for every $i\in L_3$, and $\sum_{i\in L_3}\lambda_i \ge |L_1|$, yielding that 
\begin{align}
    \label{eq:lambda}
    \sum_{i\in L_3}(\lambda_i-1) \ge \lambda.
\end{align}
{Recall that $G[tZ_i\setminus Z_1^i]$ is $K_{1, t}$-free by Claim~\ref{clm:clean0} (2). Since $Z_i^i=Z_i$ is an independent set, it follows that 
$e_G(\tZ_i\setminus Z_1^i)\le (t-1)|\tZ_i\setminus (Z_i\cup Z_1^i)|$.
Together with~\eqref{eq:toZii}, this gives
\[
e_G(\tZ_i\setminus Z_1^i)+e_G(Z_1^i, Z_i)\le (t-1)|\tZ_i\setminus (Z_i\cup Z_1^i)|+(t-1)|Z_1^i|=(t-1)(z_i-n).
\]
Therefore, for $i\in L_3$, writing $\rho_i:=e_T(Z_1^i,\tZ_i\setminus (Z_i\cup Z_1^i))$, we have
\[
e_G(\tZ_i)=e_G(\tZ_i\setminus Z_1^i)+e_G(Z_1^i, Z_i)+e_G(Z_1^i,\tZ_i\setminus (Z_i\cup Z_1^i)\le (t-1)(z_i-n)+\rho_i.
\]
Moreover, the definition of $\alpha_i$ implies that 
\[
\begin{aligned}
e_G(W_i') + e_G(\tZ_i\setminus Z_i, W_i')  - \alpha_i &\le e_T(W_i') + |\tZ_i\setminus Z_i| |W_i'|  - \alpha_i = - e_T(\tZ_i\setminus Z_i) \\
&= - \rho_i - e_T(\tZ_i\setminus( Z_i\cup Z_1)) \\
&\le -\rho_i -\binom{\lambda_i}{2}\le -\rho_i+ 1-\lambda_i.
\end{aligned}
\]
Finally, by~\eqref{eq:toZii}, we have $e_G(Z_i, W_i')\le (t-1)w_i$ for $i\in L_3$.}
Combining all these inequalities together, we obtain that, for $i\in L_3$,
\[
\begin{aligned}
e_G(\tU_i)-\alpha_i &= e_G(\tZ_i) + e_G(Z_i, W_i') + e_G(W_i') + e_G(\tZ_i\setminus Z_i, W_i')  - \alpha_i \\
&\le (t-1)(z_i-n+w_i) + (1-\lambda_i).
\end{aligned}
\]
It follows that $\sum_{i\in L_3}(e_G(\tU_i)-\alpha_i)\le \sum_{i\in L_3}(t-1)(z_i-n+w_i) -\lambda$ by using \eqref{eq:lambda}.

Using $\sum_{i=2}^r (z_i-n) = n-z_1$ and $\sum_{i=2}^{r}w_i= bn-w_1$, we derive that
\[
\begin{aligned}
 \sum_{i\in L_1} \min\{z_i, w_i\} + \sum_{i\in L_2\cup L_3} (z_i-n+w_i)
=& \sum_{i\in L_1} \min\{n-w_i, n-z_i\} + \sum_{i=2}^r (z_i-n+w_i)  \\
=& \sum_{i\in L_1} \min\{n-w_i, n-z_i\}+bn-w_1 + n-z_1.
\end{aligned}
\]
Therefore, 
\[
\begin{aligned}
\sum_{i=2}^{r} (e_G(\tU_i)-\alpha_i) &\le \sum_{i\in L_1} (t-1)\min\{z_i, w_i\} + \sum_{i\in L_2\cup L_3} (t-1)(z_i-n+w_i) -\lambda \\
&= (t-1)(bn-w_1 + n-z_1) + \sum_{i\in L_1} (t-1)\min\{n-w_i, n-z_i\} -\lambda.
\end{aligned}
\]

Finally, we work on the $\beta_i$'s and recall that $\beta_i=\sum_{j\in L\setminus \{i\}}|Z_j^i| (|\tZ_j\setminus Z_j|+|W_j'|-n + |Z_i\setminus \tZ_i|)$.
For $i\in \{1\}\cup L_1\cup L_2$, as $\tZ_i\setminus Z_i = \emptyset$, $\beta_i = 0$. For $i\in L_3$, as $|W_i|< 2\eps n$, we have $Z_i\setminus \tZ_i = \emptyset$; for any $j\in L\setminus \{i\}$, we have $\tZ_j\setminus Z_j = \emptyset$ and $W'_j\subseteq V_{q_j}$ again by Claim~\ref{eq:Wiepsn}. Hence $\beta_i=\sum_{j\in L\setminus \{i\}}|Z_j^i| (|W_j'|-n )\le 0$ because $|W'_j|\le n$. 
%
It then follows that (noting that $L\cap L_3 = \emptyset$)
\begin{align*}
    \sum_{i\ge 1}\beta_i &= \sum_{i\in L_3}\beta_i = \sum_{i\in L_3}
\sum_{j\in L\setminus\{i\}}|Z_j^i| (|W_j'| -n) \\
&= \sum_{j\in L}\sum_{i\in L_3\setminus \{j\}} |Z_j^i| (|W_j'| -n) 
= \sum_{j\in L}(n-z_j) (w_j -n).
\end{align*}
Note that $1\in L$ by \eqref{eq:W1} and $(n-z_1) (w_1 -n)\le 0$ by Claim~\ref{clm:clean}. Furthermore, since $n-z_j=0$ for $j\in L_2$, 
it follows that $\sum_{i\ge 1}\beta_i = \sum_{j\in L_1}(n-z_j) (w_j -n)$.

Recall that $e_G(\tU_1)= f_1\le z_t(z_1, w_1)$.
{By \ref{item:zgap}, we have $z_t(z_1, w_1) + (t-1)(n-z_1 + n-w_1)\le z_t(n)$}.
Thus, combining these estimates together, by~\eqref{eq:eG2}, we get
\begin{align}\label{eq:end}
    e(G)\le t_r(k) n^2 + \sum_{i\in [r]}(e_G(\tU_i)-\alpha_i+\beta_i)\le t_r(k) n^2 + z_t(n) + (t-1)(b-1)n + y - \lambda,
\end{align}
where $y:=  \sum_{i\in L_1} ((t-1)\min\{n-w_i, n-z_i\} - (n-z_i)(n-w_i))$.
For each $i\in L_1$, let $y_i:=\min\{n-w_i, n-z_i\}$ and $y_i':=\max\{n-w_i, n-z_i\}$. 
Then $y_i\le y_i'$ and thus,
\[
(t-1)\min\{n-w_i, n-z_i\} - (n-z_i)(n-w_i) = y_i (t-1- y_i') \le \lfloor (t-1)^2/4 \rfloor \le 1.
\]
{Since $L_1\subseteq L\setminus \{1\}$, we have $|L_1|\le b-1$. 
Moreover, by Claim~\ref{eq:Wiepsn}, we have $w_i\le |W_i|+|W_i'\setminus W_i|\le  2\eps n+\gamma n\le 3\eps n$ for each $i\in L_3$.
If $|L_1\cup L_2|\le b-2$, then  
\[
bn= \sum_{i\in [r]} |W_i|\le n + (b-2)n+(r-b+1)\cdot 3\eps n < bn,
\]
a contradiction.
This implies $|L_1\cup L_2|\ge b-1$, and $|L_3|\le r-b$.
Since $|L_1| - \lambda = \min\{|L_1|, |L_3|\} \le |L_3|\le r-b$, it follows that $|L_1| - \lambda \le \min\{b-1, r-b\}$.
Consequently, as $\lfloor (t-1)^2/4 \rfloor \le 1$, we get
\[
y-\lambda\le |L_1| \lfloor (t-1)^2/4 \rfloor - \lambda \le \min\{b-1, r-b\}\lfloor (t-1)^2/4 \rfloor,
\]
}
Together with \eqref{eq:end}, it gives the desired bound $e(G)\le t_r(k) n^2 + z_t(n) + (t-1)(b-1)n + \min\{b-1, r-b\}\lfloor (t-1)^2/4 \rfloor=g(n,r,k,t)$. 
This completes the proof of Theorem~\ref{thm:lb} for $k<2r$.
\end{proof}



\section*{Acknowledgments}
The authors would like to thank Junxue Zhang for creating a figure and helpful comments.


\bibliographystyle{abbrv}
\bibliography{Jan2024}


\end{document}